\documentclass[a4paper,12pt]{amsart}
\usepackage{amssymb}
\usepackage{mathrsfs}
\usepackage{color}
\usepackage{epsfig}
\usepackage{graphicx}
\usepackage[all]{xy}

\theoremstyle{plain}
\newtheorem{theorem}{Theorem}[section]
\newtheorem{prop}[theorem]{Proposition}

\newtheorem{lemma}[theorem]{Lemma}

\theoremstyle{definition}
\newtheorem{remark}[theorem]{Remark}

\newtheorem{definition}[theorem]{Definition}
\newtheorem{example}[theorem]{Example}

\newtheorem{assump}[theorem]{Assumption}

\newtheorem{assumption}[theorem]{Assumption}

\newcommand{\C}{\mathbb{C}}
\newcommand{\R}{\mathbb{R}}
\newcommand{\Z}{\mathbb{Z}}

\newcommand{\CP}{\mathbb{C}\mathrm{P}}

\newcommand{\g}{\mathfrak g}
\newcommand{\Tt}{\mathfrak t}
\newcommand{\h}{\mathfrak h}
\newcommand{\SF}{\mathscr F}

\newcommand{\vep}{\varepsilon}
\newcommand{\ind}{{\rm ind}}

\newcommand{\U}{\mathcal U}

\renewcommand{\tilde}{\widetilde}
\renewcommand{\setminus}{\smallsetminus}
\newcommand{\nin}{/\kern-2.1ex\in}

\def\<{\left\langle}
\def\>{\right\rangle}
\def\End{\operatorname{End}}

\def\ind{\operatorname{ind}}

\numberwithin{equation}{section}

\title[A Danilov-type formula for toric origami manifolds]{A Danilov-type formula for toric origami manifolds \\ via localization of index}
\date{}
\author{Hajime Fujita}

\thanks{$^1$Partly supported by Grant-in-Aid for Young Scientists (B) 26800045.}

\address{Department of Mathematical and Physical Sciences Japan Women's University 2-8-1 Mejirodai, Bunkyo-ku Tokyo, 112-8681 Japan}
\email{fujitah@fc.jwu.ac.jp}

\begin{document}

\maketitle

\begin{abstract}
We give a direct geometric proof of a Danilov-type formula for toric origami manifolds by using the localization of Riemann-Roch number. 
\end{abstract}

\tableofcontents
\section{Introduction}
A {\it symplectic toric manifold} is a symplectic manifold on which a half dimensional torus $T$ acts in an effective Hamiltonian way. A famous theorem of Delzant \cite{Delzant} says that there is one-to-one correspondence between the set of (compact connected) symplectic toric manifolds and the set of simple polytopes called {\it Delzant polytopes} (see \cite{Guilleminbook}) via moment maps. Therefore, several properties of symplectic toric manifolds, such as the symplectic volume and the ring structure of the (equivariant) cohomology and so on, can be detected from the Delzant polytopes. In view of the geometric quantization of symplectic manifolds we are interested in the {\it Riemann-Roch numbers}. The Riemann-Roch number $RR(M,L)$ is an invariant of a compact symplectic manifold $(M,\omega)$ with a {\it pre-quantizing line bundle} $(L,\nabla)$, a pair consisting of a Hermitian line bundle $L$ and a Hermitian connection $\nabla$ whose curvature form is equal to $-\sqrt{-1}\omega$, which is defined as follows.  We fix an $\omega$-compatible almost complex structure and then it determines a spin$^c$-structure of $M$ and we have a spin$^c$-Dirac operator $D$ with coefficients in $L$. We define an integer $RR(M,L)$ as the analytic index of the spin$^c$-Dirac operator: 
$$
RR(M,L):=\ind(D). 
$$
If a compact Lie group $G$ acts on $M$ preserving all the data, $\omega$, $(L,\nabla)$ and $D$, then the index becomes a virtual representation of $G$, an element of the character ring $R(G)$. In this case the Riemann-Roch number is called the {\it Riemann-Roch character} or the {\it $G$-equivariant Riemann-Roch number} and is denoted by $RR_G(M,L)$.
Such a procedure is called {\it spin$^c$-quantization} (\cite{Silva-Karshon-Tolman}\cite{Fuchs}\cite{Paradanspinc}) nowadays and considered as a quantization of spin$^c$-manifolds. 
When $(X,\omega)$ is a symplectic toric manifold with the action a torus $T$ we can choose an almost complex structure so that it is integrable and invariant under the action of the torus $T$. Then $L$ has a structure of a holomorphic line bundle and the Riemann-Roch number is equal to the dimension of $H^0(X,L)$, the space of holomorphic sections of $L$. Moreover when we consider a lift of the torus action to the pre-quantizing line bundle, $RR_T(X,L)=H^0(X,L)$ becomes a representation of the torus $T$.  Classical theorem of Danilov \cite{Danilov} says that the representation $RR_T(X,L)$ can be described in terms of the integral points in the Delzant polytope. Precisely we have 
\begin{equation}\label{originaldanilov}
RR_T(X,L)=\bigoplus_{\xi\in\mu(M)\cap{\Tt}_{\Z}^*}\C_{(\xi)}, 
\end{equation}
where $\mu$ is the moment map,  $\Tt_{\Z}^*$ is the integral weight lattice in the dual of the Lie algebra of $T$ and $\C_{(\xi)}$ is the representation of the torus associated with the integral weight $\xi\in{\Tt}_{\Z}^*$. Though Danilov's original proof was based on an algebraic geometric setting, a proof in the symplectic geometric setting is also known. See \cite{Hamiltontoric} for example.   

A {\it folded symplectic manifold} introduced by Cannas da Silva, Guillemin and Woodward in \cite{SilvaGuilleminWoodward} is a pair consisting of an even-dimensional smooth manifold and a closed 2-form which may degenerate in a transverse way and it is called the {\it folded symplectic form}. When the degenerate locus (which becomes a hypersurface and called the {\it fold}) has a structure of a circle bundle whose vertical tangent bundle coincides with the degenerate direction of the folded symplectic form, the folded symplectic manifold is called an {\it origami manifold}. By definition a folded symplectic manifold (resp. origami manifold) is a generalization of a symplectic manifold, and several notions and studies in symplectic geometry are generalized to the folded symplectic (resp. origami) case, such as pre-quantizing line bundle,  Hamiltonian group action, moment map, convexity property and so on. It is known that a folded symplectic manifold is not orientable in general, and hence it does not admit an almost complex structure, however, if it is orientable, then it admits a stable almost complex structure as shown in \cite[Theorem~2]{SilvaGuilleminWoodward}. Since the stable almost complex structure determines a spin$^c$-structure, we can define its spin$^c$-quantization by the index of spin$^c$-Dirac operator. If the folded symplectic manifold is equipped with a Hamiltonian group action, then it becomes a virtual representation and is also called the Riemann-Roch character.  In particular the spin$^c$-quantization of a toric origami manifold is a virtual representation of the torus. 

In this paper we give a proof of the following generalization of Danilov's formula (\ref{originaldanilov}) for spin$^c$-quantization of toric origami manifolds by making use of the localization theorem of index developed in \cite{Fujita-Furuta-Yoshida1, Fujita-Furuta-Yoshida2}. 

\medskip

\noindent
{\bf Theorem} (Theorem~\ref{origamiDanilov}). \ 
{\it Let $(M,\omega)$ be an oriented toric origami manifold with the action of a torus $T$ and a $T$-equivariant pre-quantizing line bundle $(L,\nabla)$. Then we have 
$$
RR_T(M,L)=\bigoplus_{\xi\in\mu(M^+)\cap{\Tt}_{\Z}^*}\C_{(\xi)}-\bigoplus_{\xi\in\mu(M^-)\cap{\Tt}_{\Z}^*}\C_{(\xi)}
$$as elements in the character ring of $T$. }

\medskip

\noindent
Precise statement and notations are explained in the subsequent sections. 
The formula itself can be obtained as a consequence of the cobordism theorem \cite[Theorem~4.1]{Silva-Guillemin-Pires} and Danilov's formula (\ref{originaldanilov}) for symplectic toric manifolds. There is an another possible approach which uses the theory of {\it multi-fans} introduced by Hattori and Masuda \cite{Hattori-Masuda1}. Masuda and Park showed in \cite{Masuda-Park} that one can associate a multi-fan for each oriented toric origami manifold. In view of the theory of multi-fans the above formula can be considered as a special case of the equivariant index formula \cite[Theorem~11.1]{Hattori-Masuda1}, which is based on the fixed point formula. In contrast to these proofs, our proof is direct and geometric, which detects the contribution of each lattice point directly.  Once we construct a geometric structure which we call an {\it acyclic compatible system} on an open subset of the manifold, then the index of Dirac operator is localized at the complement of the open subset by the localization formula in \cite{Fujita-Furuta-Yoshida2}. In this paper we construct an acyclic compatible system on the complement of the inverse image of the lattice points and the fold for toric origami manifolds. It implies that the Riemann-Roch character is equal to the sum of contributions of the lattice points and the fold. We show that the contribution of the lattice point $\xi$ is equal to $\C_{(\xi)}$ with sign determined by the orientation and the contribution of the fold is zero. Our proof does not rely on neither the original Danilov's formula nor the fixed point formula. In fact, as a special case, our proof gives a new direct proof of Danilov's formula for symplectic toric manifolds. Note that there is an another generalization of the formula (\ref{originaldanilov}) by Karshon and Tolman \cite{KarshonTolman}. They gave a formula for toric manifolds with a torus invariant {\it presymplectic form}. Though their proof is based on the holomorphic structure of toric manifolds, our proof does not use such rigid structure and it is topological and flexible. 

This paper is organized as follows. In Section~\ref{Folded symplectic forms and  toric origami manifolds} we summarize several known facts about folded symplectic manifolds, origami manifolds and toric origami manifolds, which we use in this paper. The convexity theorem for toric origami manifolds (Theorem~\ref{origamiconvexity}) is essential for us. 
In Section~\ref{Stable almost complex structure and Clifford module bundle} we discuss stable almost complex structures on folded symplectic manifolds. We construct a $\Z/2$-graded Clifford module bundle in terms of the stable almost complex structure. 
In Section~\ref{Compatible fibration on toric origami manifolds} we construct a structure of {\it (good) compatible fibration} on toric origami manifolds,  which is a family of torus fibrations (foliations) with specific compatibility condition introduced in \cite{Fujita-Furuta-Yoshida2}.  The construction is based on an open covering of the convex polytope associated with the natural stratification of the polytope with respect to the dimension of the faces. Strictly speaking there exist {\it cracks} on which we can not extend the compatible fibration keeping the compatibility condition. Though the crack causes an extra contribution to the Riemann-Roch character, we show that it is equal to 0. In Section~\ref{Compatible system on toric origami manifolds} we construct a {\it compatible system} on the compatible fibration of the toric origami manifolds, which is a family of Dirac-type operators along the fibers of the compatible fibration with specific anti-commutativity. In \cite{Fujita-Furuta-Yoshida2} the authors had already constructed compatible system for Hamiltonian torus manifolds, and our construction for the complement of the fold is based on that. On the other hand a neighbourhood of the fold has a structure of a quotient of the product of the fold and the cylinder with the standard folded symplectic structure by a natural $S^1$-action. We use this structure to define the Dirac-type operator along fibers near the fold. To discuss the localization it is essential to investigate the {\it acyclicity} of the compatible system. The fundamental property of the moment map says that it is acyclic outside the lattice points and the fold. In Section~\ref{Localization formula and Danilov type formula} we explain the localization formula of the Riemann-Roch character by making use of the acyclic compatible system. 
In Section~\ref{Computation of the local contribution} we compute the local contribution of the crack, lattice points and the fold. We first consider the symplectic toric case, i.e., origami manifolds with the empty fold, and compute the local contribution. We use a decomposition of a neighbourhood of the fiber, the inverse image of the lattice point, into the product of the cotangent bundle of the fiber and the normal direction of the symplectic submanifold containing the fiber. We apply the product formula (\cite[Theorem~8.8]{Fujita-Furuta-Yoshida2}) to the neighbourhood of the fiber. 
The vanishing of the contribution from the fold follows from the product structure of a neighbourhood of the fold. 
The last three sections are appendixes. In Appendix~\ref{Acyclic compatible systems and their local indices} we give a brief summary of the theory of local index following \cite{Fujita-Furuta-Yoshida2, Fujita-Furuta-Yoshida3} and \cite{Fujitacobinv}. In Appendix~\ref{A formula of local indeices of vector spaces} we show a useful formula of local indices of vector spaces, which will be essential in the proof of Lemma~\ref{locinddisc} and Lemma~\ref{locindcylinder}. In Appendix~\ref{A computation of local index of the folded cylinder} we give a direct computation of the local index of the folded cylinder and show that it is equal to $0$. We use this result to show that vanishing of the contribution from the fold. 


\section{Folded symplectic forms and  toric origami manifolds}
\label{Folded symplectic forms and  toric origami manifolds}
\subsection{Folded symplectic forms and  origami manifolds}
In this section we recall basic definitions and facts on folded symplectic manifolds and origami manifolds. 
Details can be found in \cite{Silva-Guillemin-Pires}, \cite{SilvaGuilleminWoodward}, \cite{Holm-Pires} and \cite{Masuda-Park}. 

A folded symplectic form $\omega$ on a smooth $2n$-dimensional manifold $M$ is a closed 2-form whose top power $\omega^n$ vanishes transversally on a submanifold $Z$ and whose restriction to $Z$ has maximal rank. In this case $Z$ is a hypersurface in $M$ and is called the {\it folding hypersurface} or {\it fold}. 
The pair $(M,\omega)$ is called a {\it folded symplectic manifold} and the $2$-form $\omega$ is called a {\it folded symplectic form}. 
Let $i_Z:Z\hookrightarrow M$ be the inclusion of $Z$ into $M$. The restriction $i_Z^*\omega$ determines a line field on $Z$, called the {\it null foliation}, whose fiber at $z\in Z$ is $\ker(i_Z^*\omega_z)$. 

Suppose that $(M,\omega)$ is an oriented folded symplectic manifold with non-empty fold $Z$. Then $M\setminus Z$ is not connected and has a decomposition $M\setminus Z=M_+\sqcup M_-$, where $M_+$ (resp. $M_-$) is the union of connected components such that $\omega^n|_{M_+}$ agrees (resp. disagrees) with the given orientation of $M$. 

\begin{definition}
A folded symplectic manifold $(M,\omega)$ is called  an {\it origami manifold} if the null foliation $\ker(i_Z^*\omega)$ is the vertical tangent bundle of a principal $S^1$-bundle structure 
$\pi:Z\to B$ over $Z$ with a compact base $B$. 
\end{definition}
Note that since $B$ is compact the total space $Z$ is also compact. As in the symplectic reduction procedure, there is the unique symplectic form $\omega_B$ on $B$ satisfying $\pi^*\omega_B=i_Z^*\omega$. 
An analogue of Darboux's theorem for folded symplectic forms says that near any point $p\in Z$ there exists a coordinate chart centered at $p$ where the folded symplectic form $\omega$ can be written as  
$$
x_1dx_1\wedge dy_1+dx_2\wedge dy_2+\cdots + dx_n\wedge dy_n. 
$$
In this local description, the fold $Z$ is given by the equation $x_1=0$ and the null foliation is the line field spanned by $\frac{\partial}{\partial y_1}$.  
This local description has a global variant. 

\begin{theorem}[Theorem~1 in \cite{SilvaGuilleminWoodward} ]\label{Morsermodel} 
Let $(M,\omega)$ be an oriented origami manifold with fold $Z\to B$. 
Fix a connection 1-form $\alpha$ of $Z\to B$. Then there exists a neighbourhood ${\mathcal U}$ of $Z$ and an orientation preserving diffeomorphism $\varphi:Z\times(-\varepsilon,\varepsilon)\to {\mathcal U}$ such that 
$$
\varphi\circ\iota_0=\iota_Z
$$and 
$$
\varphi^*\omega=p_Z^*\iota^*_Z\omega+d(t^2p_Z^*\alpha), 
$$where $\iota_0: Z \to Z\times(-\varepsilon,\varepsilon)$ is the inclusion 
$z\mapsto (z,0)$ and $p_Z:Z\times(-\varepsilon,\varepsilon)\to Z$ is the natural projection. 
\end{theorem}

\begin{example}\label{exsphere2}
For a positive integer $n$ let $S^{2n}$ be the unit sphere in $\R^{2n}\oplus\R=\C^{n}\oplus \R$ with coordinates $x_1, y_1, \cdots, x_n, y_n, h$. Let $\omega$ be the restriction to $S^{2n}$ of the 2-form  $dx_1\wedge dy_1+\cdots+dx_n\wedge dy_n$ on $\R^{2n}\oplus \R$. Then $\omega$ is a folded symplectic form on $S^{2n}$ with the fold $S^{2n-1}$, the equator sphere given by $h=0$. The Hopf fibration $S^1\hookrightarrow S^{2n-1}\to \CP^{n-1}$ gives a structure of origami manifold on $(S^{2n}, \omega)$. 
\end{example}

\subsection{Hamiltonian torus actions and toric origami manifolds}

The action of a compact Lie group $G$ on an origami manifold $(M,\omega)$ is called {\it Hamiltonian} if it admits a moment map $\mu$, that is, a map $\mu:M\to {\mathfrak g}^*={\rm Lie}(G)^*$ satisfying the conditions : 
\\
$\bullet$ $\mu$ is equivariant with respect to the given action of $G$ on $M$ and the coadjoint action of $G$ on ${\mathfrak g}^*$. 
\\
\noindent 
$\bullet$ for any $v\in {\mathfrak g}$ we have $d\langle \mu, v \rangle = \iota({v^M})\omega$, where  $\langle \cdot, \cdot \rangle$ is the pairing between ${\mathfrak g}^*$ and ${\mathfrak g}$ and $\iota({v^M})\omega$ is the contraction of $\omega$ by  the induced fundamental vector field $v^M$. 

\begin{definition}
A Hamiltonian torus origami manifold $(M,\omega, T,\mu)$ (or $M$ for short) is a connected origami manifold $(M,\omega)$ equipped with an effective Hamiltonian action of a torus $T$ with a choice of a corresponding moment map $\mu$. If the dimension of the torus $T$ is half of that of $M$, then we call $(M,\omega, T,\mu)$ a {\it toric origami manifold}. 
\end{definition}

If the fold $Z$ is empty, a Hamiltonian torus origami manifold is a Hamiltonian  torus manifold in the usual sense. 
The following is an origami analogue of the famous convexity theorem for Hamiltonian torus manifolds. 

\begin{theorem}[Theorem~3.2 in \cite{Silva-Guillemin-Pires}]\label{origamiconvexity}
Let $(M,\omega,T,\mu)$ be  a connected compact origami manifold with null fibration $\pi:Z\to B$ and a Hamiltonian torus action of a torus $T$ with moment map $\mu$. Then : \\
{\rm (a)} The image $\mu(M)$ is the union of a finite number of convex polytopes $\Delta_1,\cdots, \Delta_N$ in the dual of the Lie algebra ${\mathfrak t}^*$, each of which is the image of the moment map restricted to the closure of a connected component of $M\setminus Z$. 
\\
{\rm (b)} Over each connected component $Z'$ of $Z$, the null fibration is given by a subgroup of $T$ if and only if $\mu(Z')$ is a facet of each of the one or two polytopes corresponding to the neighbourhood(s) of $M\setminus Z$, and when those are two polytopes $\Delta_1$ and $\Delta_2$ there exists an open subset $\tilde \Delta_{Z'}$ containing $\mu(Z')$ such that ${\tilde \Delta_{Z'}}\cap\Delta_1={\tilde \Delta_{Z'}}\cap \Delta_2$. 

We call such images $\mu(M)$ origami polytopes. 
\end{theorem}

\begin{example}\label{exsphere3}
Consider the origami manifold $(S^{2n}, \omega)$ given in Example~\ref{exsphere2}. Let $T:=(S^1)^n$ be the $n$-dimensional torus. Then the action of $T$ on $S^{2n}$ given by 
$$ 
(t_1,\ldots, t_n)\cdot (z_1, \ldots, z_1, h):=(t_1z_1,\ldots, t_nz_n, h)
$$for $(t_1,\ldots, t_n)\in T$ and $(z_1, \ldots, z_1,h)\in S^{2n}\subset \C^n\oplus \R$ is Hamiltonian (in fact, toric) action with the moment map $\mu:S^{2n}\to \R^n$, 
$$
\mu(z_1,\ldots,z_n,h):=\frac{1}{2}(|z_1|^2, \ldots, |z_n|^2). 
$$The image of $\mu$ is the union of two copies of the $n$-simplex, $\xi_1,\cdots,\xi_n\geq 0$, $\xi_1+\cdots+\xi_n\leq 1/2$,  and the image of fold $S^{2n-1}$ is the \lq\lq hypotenuse\rq\rq, $\xi_1+\cdots +\xi_n=1/2$. See Figure~\ref{pics4} for the case of $n=2$. 

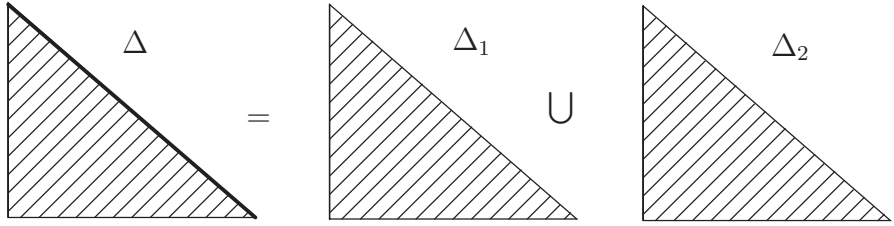
\begin{figure}[h]
{\unitlength 0.1in
\begin{picture}( 45.6800, 11.3600)( 16.0000,-27.3600)
%
{\color[named]{Black}{%
\special{pn 8}%
\special{pa 1600 1600}%
\special{pa 1600 2720}%
\special{pa 2880 2720}%
\special{pa 1600 1600}%
\special{pa 1600 2720}%
\special{fp}%
}}%
\put(28.3200,-22.4000){\makebox(0,0)[lb]{$=$}}%
%
{\color[named]{Black}{%
\special{pn 8}%
\special{pa 2232 2152}%
\special{pa 1664 2720}%
\special{fp}%
\special{pa 2280 2200}%
\special{pa 1760 2720}%
\special{fp}%
\special{pa 2336 2240}%
\special{pa 1856 2720}%
\special{fp}%
\special{pa 2384 2288}%
\special{pa 1952 2720}%
\special{fp}%
\special{pa 2432 2336}%
\special{pa 2048 2720}%
\special{fp}%
\special{pa 2488 2376}%
\special{pa 2144 2720}%
\special{fp}%
\special{pa 2536 2424}%
\special{pa 2240 2720}%
\special{fp}%
\special{pa 2592 2464}%
\special{pa 2336 2720}%
\special{fp}%
\special{pa 2640 2512}%
\special{pa 2432 2720}%
\special{fp}%
\special{pa 2688 2560}%
\special{pa 2528 2720}%
\special{fp}%
\special{pa 2744 2600}%
\special{pa 2624 2720}%
\special{fp}%
\special{pa 2792 2648}%
\special{pa 2720 2720}%
\special{fp}%
\special{pa 2848 2688}%
\special{pa 2816 2720}%
\special{fp}%
\special{pa 2176 2112}%
\special{pa 1600 2688}%
\special{fp}%
\special{pa 2128 2064}%
\special{pa 1600 2592}%
\special{fp}%
\special{pa 2080 2016}%
\special{pa 1600 2496}%
\special{fp}%
\special{pa 2024 1976}%
\special{pa 1600 2400}%
\special{fp}%
\special{pa 1976 1928}%
\special{pa 1600 2304}%
\special{fp}%
\special{pa 1920 1888}%
\special{pa 1600 2208}%
\special{fp}%
\special{pa 1872 1840}%
\special{pa 1600 2112}%
\special{fp}%
\special{pa 1824 1792}%
\special{pa 1600 2016}%
\special{fp}%
\special{pa 1768 1752}%
\special{pa 1600 1920}%
\special{fp}%
\special{pa 1720 1704}%
\special{pa 1600 1824}%
\special{fp}%
\special{pa 1664 1664}%
\special{pa 1600 1728}%
\special{fp}%
}}%
%
{\color[named]{Black}{%
\special{pn 8}%
\special{pa 3264 1608}%
\special{pa 3264 2728}%
\special{pa 4544 2728}%
\special{pa 3264 1608}%
\special{pa 3264 2728}%
\special{fp}%
}}%
%
{\color[named]{Black}{%
\special{pn 8}%
\special{pa 3896 2160}%
\special{pa 3328 2728}%
\special{fp}%
\special{pa 3944 2208}%
\special{pa 3424 2728}%
\special{fp}%
\special{pa 4000 2248}%
\special{pa 3520 2728}%
\special{fp}%
\special{pa 4048 2296}%
\special{pa 3616 2728}%
\special{fp}%
\special{pa 4096 2344}%
\special{pa 3712 2728}%
\special{fp}%
\special{pa 4152 2384}%
\special{pa 3808 2728}%
\special{fp}%
\special{pa 4200 2432}%
\special{pa 3904 2728}%
\special{fp}%
\special{pa 4256 2472}%
\special{pa 4000 2728}%
\special{fp}%
\special{pa 4304 2520}%
\special{pa 4096 2728}%
\special{fp}%
\special{pa 4352 2568}%
\special{pa 4192 2728}%
\special{fp}%
\special{pa 4408 2608}%
\special{pa 4288 2728}%
\special{fp}%
\special{pa 4456 2656}%
\special{pa 4384 2728}%
\special{fp}%
\special{pa 4512 2696}%
\special{pa 4480 2728}%
\special{fp}%
\special{pa 3840 2120}%
\special{pa 3264 2696}%
\special{fp}%
\special{pa 3792 2072}%
\special{pa 3264 2600}%
\special{fp}%
\special{pa 3744 2024}%
\special{pa 3264 2504}%
\special{fp}%
\special{pa 3688 1984}%
\special{pa 3264 2408}%
\special{fp}%
\special{pa 3640 1936}%
\special{pa 3264 2312}%
\special{fp}%
\special{pa 3584 1896}%
\special{pa 3264 2216}%
\special{fp}%
\special{pa 3536 1848}%
\special{pa 3264 2120}%
\special{fp}%
\special{pa 3488 1800}%
\special{pa 3264 2024}%
\special{fp}%
\special{pa 3432 1760}%
\special{pa 3264 1928}%
\special{fp}%
\special{pa 3384 1712}%
\special{pa 3264 1832}%
\special{fp}%
\special{pa 3328 1672}%
\special{pa 3264 1736}%
\special{fp}%
}}%
%
{\color[named]{Black}{%
\special{pn 8}%
\special{pa 4888 1616}%
\special{pa 4888 2736}%
\special{pa 6168 2736}%
\special{pa 4888 1616}%
\special{pa 4888 2736}%
\special{fp}%
}}%
%
{\color[named]{Black}{%
\special{pn 8}%
\special{pa 5520 2168}%
\special{pa 4952 2736}%
\special{fp}%
\special{pa 5568 2216}%
\special{pa 5048 2736}%
\special{fp}%
\special{pa 5624 2256}%
\special{pa 5144 2736}%
\special{fp}%
\special{pa 5672 2304}%
\special{pa 5240 2736}%
\special{fp}%
\special{pa 5720 2352}%
\special{pa 5336 2736}%
\special{fp}%
\special{pa 5776 2392}%
\special{pa 5432 2736}%
\special{fp}%
\special{pa 5824 2440}%
\special{pa 5528 2736}%
\special{fp}%
\special{pa 5880 2480}%
\special{pa 5624 2736}%
\special{fp}%
\special{pa 5928 2528}%
\special{pa 5720 2736}%
\special{fp}%
\special{pa 5976 2576}%
\special{pa 5816 2736}%
\special{fp}%
\special{pa 6032 2616}%
\special{pa 5912 2736}%
\special{fp}%
\special{pa 6080 2664}%
\special{pa 6008 2736}%
\special{fp}%
\special{pa 6136 2704}%
\special{pa 6104 2736}%
\special{fp}%
\special{pa 5464 2128}%
\special{pa 4888 2704}%
\special{fp}%
\special{pa 5416 2080}%
\special{pa 4888 2608}%
\special{fp}%
\special{pa 5368 2032}%
\special{pa 4888 2512}%
\special{fp}%
\special{pa 5312 1992}%
\special{pa 4888 2416}%
\special{fp}%
\special{pa 5264 1944}%
\special{pa 4888 2320}%
\special{fp}%
\special{pa 5208 1904}%
\special{pa 4888 2224}%
\special{fp}%
\special{pa 5160 1856}%
\special{pa 4888 2128}%
\special{fp}%
\special{pa 5112 1808}%
\special{pa 4888 2032}%
\special{fp}%
\special{pa 5056 1768}%
\special{pa 4888 1936}%
\special{fp}%
\special{pa 5008 1720}%
\special{pa 4888 1840}%
\special{fp}%
\special{pa 4952 1680}%
\special{pa 4888 1744}%
\special{fp}%
}}%
\put(43.9200,-22.3200){\makebox(0,0)[lb]{$\bigcup$}}%
%
{\color[named]{Black}{%
\special{pn 20}%
\special{pa 1600 1608}%
\special{pa 2880 2720}%
\special{fp}%
}}%
\put(21.9000,-18.6000){\makebox(0,0)[lb]{$\Delta$}}%
\put(39.0000,-18.7000){\makebox(0,0)[lb]{$\Delta_1$}}%
\put(55.5000,-19.0000){\makebox(0,0)[lb]{$\Delta_2$}}%
\end{picture}}%
\caption{An origami polytope for $S^4$}\label{pics4}
\end{figure}

\end{example}

\section{Stable almost complex structure and Clifford module bundle}
\label{Stable almost complex structure and Clifford module bundle}
Let $(M,\omega)$ be a $2n$-dimensional oriented folded symplectic manifold with fold $Z$ and ${\mathcal U}$ an open neighbourhood of $Z$ as in Theorem~\ref{Morsermodel}. 
Let $M_+$ (resp. $M_-$) be the union of connected components of $M\setminus Z$ such that $\omega^n|_{M_+}$ agrees (resp. disagrees) with the given orientation of $M$. 
In \cite{SilvaGuilleminWoodward}, it was shown that $M$ has a stable almost complex structure. More precisely the following holds. 
\begin{theorem}[Theorem~2 in \cite{SilvaGuilleminWoodward}]\label{stablealmostcomplex}
There exixts an almost complex structure $\tilde J$ on the real $(2n+2)$-dimensional vector bundle $TM\oplus\R^2$ , and a $\C$-linear isomorphism 
$$
(TM\oplus\R^2)|_{M\setminus\U}\cong T(M\setminus \U)\oplus \C. 
$$Moreover, $TM\oplus\R^2$ has a symplectic structure $\tilde\omega$ which is canonical up to homotopy, and the homotopy class of $\tilde J$ is unique provided $\tilde J$ is compatible with the natural symplectic structure on $TM\oplus \R^2$. 
\end{theorem}

\begin{remark}\label{stablecomplexstrrem}
One can see in the proof of \cite[Theorem~2]{SilvaGuilleminWoodward} that the above $\tilde J$ has the following properties. 
\begin{enumerate}
\item Let $J$ be an almost complex structure on $M\setminus Z$ which is compatible with $\omega|_{M\setminus Z}$. Then one can construct $\tilde J$ so that the following equality holds. 
\begin{equation}\label{stablestd}
\begin{cases}
 \ \tilde J|_{M_+\setminus {\mathcal U}}=J|_{M_+\setminus {\mathcal U}}\oplus (\sqrt{-1}) \\ 
 \ \tilde J|_{M_-\setminus {\mathcal U}}=J|_{M_-\setminus {\mathcal U}}\oplus (-\sqrt{-1}). 
\end{cases}
\end{equation}
\item By using a connection of the principal $S^1$-bundle $Z\to B$ we have the splitting of the tangent bundle $TZ\cong \pi^*TB\oplus T_{\pi}Z$, where $T_{\pi}Z$ is the tangent bundle along the fiber, which is a real line bundle over $Z$. Since $T\U$ is oriented, and hence, $TZ$ is also oriented, the fact that $B$ is a symplectic manifold implies that $T_{\pi}Z$ is an orientable.  In particular, $T_{\pi}Z$ is trivial real line bundle. Under these identifications the almost complex structure $\tilde J|_{\U}$ in Theorem~\ref{stablealmostcomplex}  on $T\U\oplus \R^2\cong \pi^*TB\oplus T_{\pi}Z\oplus\R\oplus\R^2$ can be taken as the direct sum of almost complex structures on the symplectic vector bundle $\pi^*TB$ and the trivial bundle $T_{\pi}Z\oplus\R\oplus\R^2$ of real rank 4. 
\item If a compact Lie group $G$ acts on $(M,\omega)$, then we can take $\tilde J$ to be $G$-invariant. In fact we will use such an invariant $\tilde J$ in the subsequent sections. 
\end{enumerate}
\end{remark}

By using $\tilde J$ and $\tilde \omega$, we have a Riemannian metric on $TM\oplus \R^2$, and $TM$ is equipped with the metric as a subbundle of $TM\oplus\R^2$. 
Moreover the stable almost complex structure induces a spin$^c$-structure on $M$.  Now we construct a Clifford module bundle over $TM$ in terms of this stable almost complex structure.

We first explain the construction for the vector space case. 
Let $E$ be an even dimensional Euclidean vector space. 
Suppose that a complex structure $J_{\tilde E}$ on $\tilde E:=E\oplus\R^e$ which preserves the metric on $\tilde E$ is given for a non-negative (even) integer $e$. 
By using $J_{\tilde E}$ we have a $\Z/2$-graded $Cl(\tilde E)=Cl(E)\otimes Cl(\R^e)$-module $W_{\tilde E}:=\wedge^{\bullet}_{\C}\tilde E$, the exterior product algebra of the Hermitian vector space $\tilde E$. 
The Clifford action of $Cl(\tilde E)$ is defined by the wedge product and the interior product. We define $W_E$ as the set of all linear maps from an irreducible representation $W_e$ of the Clifford algebra $Cl_e:=Cl(\R^e)$ to $W_{\tilde E}$ which commute with the Clifford action of $Cl_e$, 
$$
W_E:={\rm Hom}_{Cl_e}(W_e,W_{\tilde E}), 
$$where $Cl_e$ acts on $W_{\tilde E}$ by using the inclusion $Cl_e\hookrightarrow Cl(E\oplus\R^e)$. Note that $W_E$ is equipped with the Clifford action of $Cl(E)$ by 
$$
\alpha\cdot\phi : v\mapsto \alpha \phi(v)
$$for $\alpha\in Cl(E)$ and $v\in W_e$ using the inclusion $Cl(E)\hookrightarrow Cl(E\oplus\R^e)$. 

\begin{lemma}
$W_E$ is an irreducible $\Z/2$-graded $Cl(E)$-module. 
\end{lemma}
\begin{proof}
Suppose that $E$ is equipped with an almost complex structure $J_E$ and $J_{\tilde E}$ is the direct sum of $J_E$ and the standard complex structure $\sqrt{-1}$ on $\R^e=\C^{e/2}$ (for a specific order of the basis of $\R^e$). 
In this case, one can see that $\wedge^{\bullet}_{\C}{\tilde E}=\wedge^{\bullet}_{\C}{E}\otimes\wedge_{\C}^{\bullet}\R^{e}$ and 
$$
W_E={\rm Hom}_{Cl_e}(W_e,\wedge^{\bullet}_{\C}{\tilde E})=\wedge_{\C}^{\bullet}E\otimes{\rm Hom}_{Cl_e}(W_e, \wedge_{\C}^{\bullet}\R^{e})=\wedge_{\C}^{\bullet}E.
$$
It implies that $W_E$ is an irreducible $Cl(E)$-module. 
Since any complex structure on $\tilde E$ is homotopic to the direct sum $J_E\oplus {\sqrt{-1}}$ and the irreducible representation of $Cl(E)$ is unique, we complete the proof. 
\end{proof}
By applying the above construction for an almost complex structure on $TM\oplus\R^2$ we have the $\Z/2$-graded $Cl(TM)$-module bundle 
\begin{equation}\label{Clifford}
W:={\rm Hom}_{Cl_2}(W_2, \wedge_{\C}^{\bullet}(TM\oplus\R^2)) 
\end{equation}
over $M$. Note that we have $W|_{M_{\pm}\setminus{\mathcal U}}\cong\wedge_{\C}^{\bullet}T(M_{\pm}\setminus{\mathcal U})$ by (\ref{stablestd}), which is the standard $Cl(T(M_{\pm}\setminus {\mathcal U}))$-module bundle of $M_{\pm}\setminus {\mathcal U}$. 
For any Hermitian line bundle $L$ we have an another $\Z/2$-graded $Cl(E)$-module bundle $W_L:=W\otimes L$.

\begin{definition}\label{origamiRR}
For a compact oriented origami manifold $(M,\omega)$ without boundary and a Hermitian line bundle $L$ over $M$ the {\it Riemann-Roch number} $RR(M,L)$ is defined as the index of spin$^c$-Dirac operator which acts on the smooth sections of the Clifford module bundle $W_L$: 
$$
RR(M,L):=\ind(W_L). 
$$
\end{definition}
\begin{remark}
Since any two Dirac-type operators can be joined in the space of Dirac-type operators the index $RR(M,L)=\ind(W_L)$ does not depend on the choice of the Dirac-type operators by the homotopy invariance of the analytic index. 
\end{remark}

\section{Compatible fibration on toric origami manifolds}
\label{Compatible fibration on toric origami manifolds}
In this section we construct a structure of {\it good compatible fibration} on toric origami manifolds. 
The notion of good compatible fibration is a family of torus fibrations (or more generally foliations) over an open covering of the manifold with some compatibility condition and is introduced in \cite{Fujita-Furuta-Yoshida2}.  
See also Definition~\ref{goodcompatifib}. 

\begin{assump}\label{assump}
In this section we consider a toric origami manifold $(M,\omega, T,\mu)$ satisfying the following assumptions.
\begin{itemize}
\item $M$ is connected, oriented, and compact without boundary. 
\item $(M,\omega, T,\mu)$ satisfies the condition (b) in Theorem~\ref{origamiconvexity}. Namely, the null foliation is given by a subgroup of $T$. 
\end{itemize}
\end{assump}
Suppose that $\dim M=2n$.
Let $\mu(M)=\bigcup_i\Delta_i$ be the union of convex polytopes associated with the moment map $\mu:M\to \Tt^*$. 
For each $i$ let $\Delta_{i}=\Delta_Z\cup\bigcup_{j=0}^n\bigcup_{k=1}^{m_j}\Delta_{i,k}^{(j)}$ be the stratification of $\Delta_i$, where we put\footnote{Strictly speaking we consider each connected component of $Z$. } $\Delta_Z:=\mu(Z)$ and $\{\Delta_{i,1}^{(j)},\cdots, \Delta_{i,m_j}^{(j)}\}$ is the set of all $j$-dimensional faces of $\Delta_i$ for each $j\in\{0,\cdots,n\}$. 
We take and fix a neighbourhood $\U:=Z\times (-\vep,\vep)$ of $Z$ in $M$ as in Theorem~\ref{Morsermodel} for some small $\vep>0$, and we may assume that $\overline{\U}=Z\times[-\vep,\vep,]$ and an open neighbourhood $\tilde\Delta_Z$ in Theorem~\ref{origamiconvexity}(b) has the form $\tilde \Delta_Z=\mu(\U)$. 

The construction of the good compatible fibration is divided into two parts, 
fibrations near the fold and fibrations outside the fold. 

\subsection{Torus actions near the fold}
\label{Torus actions near the fold}
We set $\U':=Z\times (-\frac{\vep}{2}, \frac{\vep}{2})$. We take $\vep>0$ small enough so that the $S^1$-aciton on $Z$ can be extended to a free $S^1$-action on $\U'$. By using this $S^1$-action we have an $S^1$-bundle structure on $\U'$ with the base space $B\times(-\frac{\vep}{2}, \frac{\vep}{2})$.

\subsection{Torus actions outside the fold}
\label{Torus actions outside the fold}
\begin{figure}[h]
\input{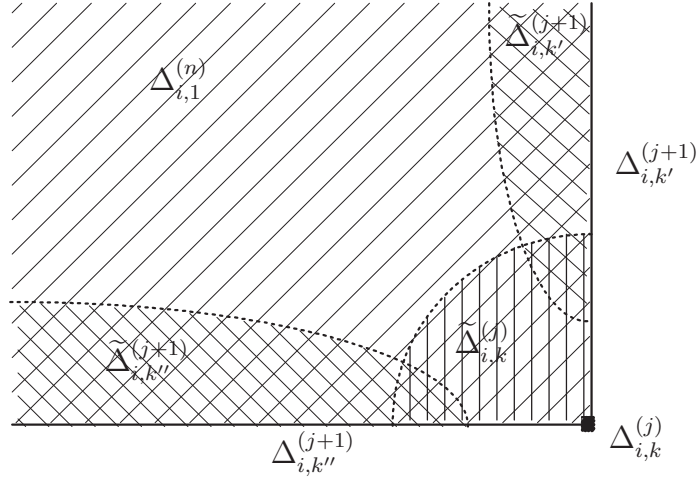}
\caption{Open covering outside the image of the fold}
\label{Open covering outside the imge of the fold}
\end{figure}

We construct a family of torus actions on $M\setminus\overline{\U}$. 
We put $\Delta_{i,Z}':=\Delta_i\setminus\mu(\overline{\U})$ and 
we first construct an open covering  
$$\Delta_{i,Z}'=\left(\bigcup_{j=0}^{n}\bigcup_{k=1}^{m_j}\tilde\Delta_{i,k}^{(j)}\right)$$ 
by the following procedure. See also Figure~\ref{Open covering outside the imge of the fold}. 

\begin{itemize}
\item[(0)] For each $k\in\{1,\cdots, m_0\}$ take a small open neighbourhood $\tilde\Delta_{i,k}^{(0)}$ of $\Delta_{i,k}^{(0)}$ in $\Delta_{i,Z}'$ so that 
$\tilde\Delta_{i,k}^{(0)}\cap\tilde\Delta_{i,k'}^{(0)}=\emptyset$ if $k\neq  k'$. 
\item[(1)] For each $k\in\{1,\cdots, m_1\}$ take a small open neighbourhood $\tilde\Delta_{i,k}^{(1)}$ of 
$$
\left(\Delta_{i,Z}'\setminus\bigcup_{k_0=1}^{m_0}\tilde\Delta_{i,k_0}^{(0)}\right)\cap \Delta_{i,k}^{(1)}
$$ in $\Delta_{i,Z}'$ so that $\tilde\Delta_{i,k}^{(1)}\cap\tilde\Delta_{i,k'}^{(1)}=\emptyset$ if $k\neq  k'$. 
\item[(2)] For each $k\in\{1,\cdots, m_2\}$ take a small open neighbourhood $\tilde\Delta_{i,k}^{(2)}$ of 
$$
\left(\Delta_{i,Z}'\setminus\bigcup_{j=0}^{1}\bigcup_{k_j=1}^{m_j}\tilde\Delta_{i,k_j}^{(j)}\right)\cap \Delta_{i,k}^{(2)}
$$ in $\Delta_{i,Z}'$ so that $\tilde\Delta_{i,k}^{(2)}\cap\tilde\Delta_{i,k'}^{(2)}=\emptyset$ if $k\neq  k'$. 
\\ 
$\vdots$
\\
\item[(n-1)] For each $k\in\{1,\cdots, m_{n-1}\}$ take a small open neighbourhood $\tilde\Delta_{i,k}^{(n-1)}$ of 
$$
\left(\Delta_{i,Z}'\setminus\bigcup_{j=0}^{n-2}\bigcup_{k_j=0}^{m_j}\tilde\Delta_{i,k_j}^{(j)}\right)\cap \Delta_{i,k}^{(n-1)}
$$ in $\Delta_{i,Z}'$ so that $\tilde\Delta_{i,k}^{(n-1)}\cap\tilde\Delta_{i,k'}^{(n-1)}=\emptyset$ if $k\neq  k'$. 
\item[(n)] We set $\tilde\Delta_{i,1}^{(n)}:={\rm int}(\Delta_i)={\rm int}(\Delta_i^{(n)})$. 
\end{itemize}

\medskip

For a toric manifold  $M\setminus Z$ it is well-known that for each $i,j,k$ there exists a subtorus $T_{i,k}^{(j)}$ of $T$ such that $\dim T_{i,k}^{(j)}=n-j$ and for any $x\in \mu^{-1}({\rm int}(\Delta_{i,k}^{(j)})) \setminus Z$ the stabilizer subgroup at $x$ is equal to $T_{i,k}^{(j)}$. 
Particularly we have $T_{i,1}^{(n)}=\{e\}$.  
We take and fix a rational metric of the Lie algebra ${\mathfrak t}$ so that for each subspace ${\mathfrak h}$ in ${\mathfrak t}$ spanned by rational vectors one can associate the orthogonal complement subgroup ${\rm exp}({\mathfrak h}^{\perp})$ as a compact subgroup of $T$. 
Let $G_{i,k}^{(j)}$ be the orthogonal complement subgroup associated with (the Lie algebra of) the stabilizer subgroup $T_{i,k}^{(j)}$.  Note that we have $G^{(n)}_{i,1}=T$. Define an open subset of $M$ by  $M_{i,k}^{(j)}:=\mu^{-1}(\tilde\Delta_{i,k}^{(j)})$, which has the natural $G_{i,k}^{(j)}$-action and the following properties. 

\begin{itemize}
\item Each $G_{i,k}^{(j)}$ acts on $M_{i,k}^{(j)}$, and all orbits of $G_{i,k}^{(j)}$-action have the maximal dimension $\dim G_{i,k}^{(j)}$.  
\item If $\tilde\Delta_{i,k}^{(j)}\cap \tilde\Delta_{i,k'}^{(j')}\neq \emptyset$, then we have $G_{i,k}^{(j)}\subset G_{i,k'}^{(j')}$ or $G_{i,k}^{(j)}\supset G_{i,k'}^{(j')}$.  
\end{itemize}

\subsection{Good compatible fibration on toric origami manifolds}
\label{Good compatible fibration on toric origami manifolds}
By taking each open subset small enough we may assume that 
$\U'\cap M_{i,k}^{(j)}=\emptyset$ for all $i,j,k$ with $j\neq n$. 
The union $\displaystyle\U'\cup\bigcup_{i,j,k}M_{i,k}^{(j)}$ is not an open covering of the whole $M$. There exist a family of compact sets, which we call the {\it crack} $C_{i,k}^Z$ defined by 
$$
C_{i,k}^{Z}:=\mu^{-1}\left(\mu(\overline\U\setminus\U')\cap\Delta_{i,k}^{(n-1)}\right). 
$$

\begin{figure}[h]
{\unitlength 0.1in%
\begin{picture}(38.8000,25.4000)(15.1000,-37.5000)%
%
\special{pn 8}%
\special{pa 2096 3670}%
\special{pa 2096 2550}%
\special{pa 3536 1430}%
\special{pa 5296 1430}%
\special{fp}%
%
\special{pn 8}%
\special{pa 2088 3590}%
\special{pa 4824 1422}%
\special{da 0.070}%
%
\special{pn 4}%
\special{pa 3600 2382}%
\special{pa 3040 1822}%
\special{fp}%
\special{pa 3688 2326}%
\special{pa 3120 1758}%
\special{fp}%
\special{pa 3768 2262}%
\special{pa 3200 1694}%
\special{fp}%
\special{pa 3848 2198}%
\special{pa 3280 1630}%
\special{fp}%
\special{pa 3928 2134}%
\special{pa 3360 1566}%
\special{fp}%
\special{pa 4008 2070}%
\special{pa 3440 1502}%
\special{fp}%
\special{pa 4088 2006}%
\special{pa 3528 1446}%
\special{fp}%
\special{pa 4168 1942}%
\special{pa 3656 1430}%
\special{fp}%
\special{pa 4248 1878}%
\special{pa 3800 1430}%
\special{fp}%
\special{pa 4328 1814}%
\special{pa 3944 1430}%
\special{fp}%
\special{pa 4408 1750}%
\special{pa 4088 1430}%
\special{fp}%
\special{pa 4488 1686}%
\special{pa 4232 1430}%
\special{fp}%
\special{pa 4568 1622}%
\special{pa 4376 1430}%
\special{fp}%
\special{pa 4648 1558}%
\special{pa 4520 1430}%
\special{fp}%
\special{pa 4728 1494}%
\special{pa 4664 1430}%
\special{fp}%
\special{pa 3520 2446}%
\special{pa 2960 1886}%
\special{fp}%
\special{pa 3440 2510}%
\special{pa 2880 1950}%
\special{fp}%
\special{pa 3360 2574}%
\special{pa 2792 2006}%
\special{fp}%
\special{pa 3280 2638}%
\special{pa 2712 2070}%
\special{fp}%
\special{pa 3200 2702}%
\special{pa 2632 2134}%
\special{fp}%
\special{pa 3120 2766}%
\special{pa 2552 2198}%
\special{fp}%
\special{pa 3040 2830}%
\special{pa 2472 2262}%
\special{fp}%
\special{pa 2960 2894}%
\special{pa 2392 2326}%
\special{fp}%
\special{pa 2880 2958}%
\special{pa 2312 2390}%
\special{fp}%
\special{pa 2800 3022}%
\special{pa 2232 2454}%
\special{fp}%
\special{pa 2720 3086}%
\special{pa 2144 2510}%
\special{fp}%
\special{pa 2640 3150}%
\special{pa 2096 2606}%
\special{fp}%
\special{pa 2560 3214}%
\special{pa 2096 2750}%
\special{fp}%
\special{pa 2480 3278}%
\special{pa 2096 2894}%
\special{fp}%
\special{pa 2400 3342}%
\special{pa 2096 3038}%
\special{fp}%
\special{pa 2320 3406}%
\special{pa 2096 3182}%
\special{fp}%
\special{pa 2240 3470}%
\special{pa 2096 3326}%
\special{fp}%
\special{pa 2160 3534}%
\special{pa 2096 3470}%
\special{fp}%
\put(26.3200,-19.0200){\makebox(0,0)[lb]{$\Delta_Z$}}%
\put(20.0000,-39.1000){\makebox(0,0)[lb]{$\Delta_{i,k}^{(n-1)}$}}%
\put(27.7000,-24.0000){\makebox(0,0)[lb]{$\mu({\mathcal U})$}}%
\put(53.9000,-15.6000){\makebox(0,0)[lb]{$\Delta_{i,k'}^{(n-1)}$}}%
%
\special{pn 20}%
\special{pa 4790 1430}%
\special{pa 4160 1430}%
\special{fp}%
\put(15.1000,-35.5000){\makebox(0,0)[lb]{$\mu(C_{i,k}^{Z})$}}%
\put(45.4000,-13.7000){\makebox(0,0)[lb]{$\mu(C_{i,k'}^{Z})$}}%
%
\special{pn 20}%
\special{pa 2090 3040}%
\special{pa 2090 3580}%
\special{fp}%
\end{picture}}%
\caption{Crack near the fold}
\label{Crack near the fold}
\end{figure}
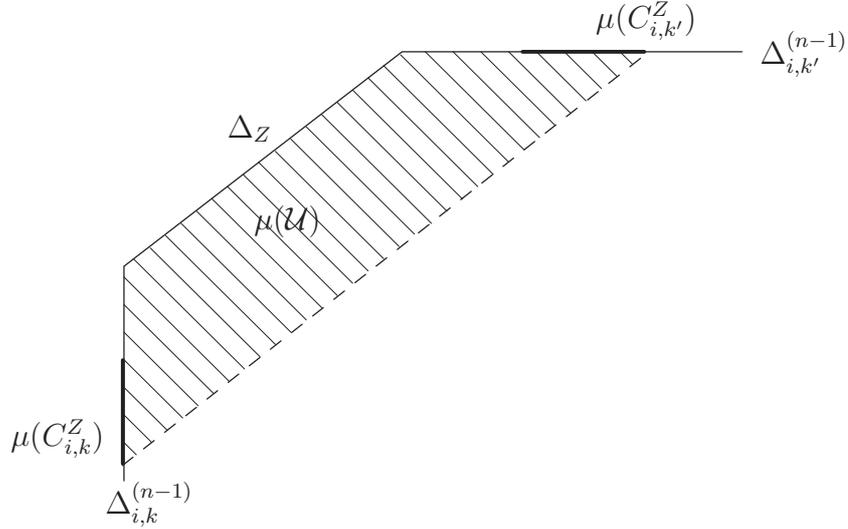

Though we do not know the way to extend the good compatible fibration across the crack, we have the following.

\begin{prop}\label{crackgoodcompatifib}
A family of open subsets $\{\U', M_{i,k}^{(j)}\}_{i,j,k}$ defines a structure of good compatible fibration (Definition~\ref{goodcompatifib}) on the complement $M\setminus \bigcup_{i,k}C_{i,k}^{Z}$. 
\end{prop}

\begin{example}\label{exshpere5}
Consider the toric origami manifold $S^4$ with the moment map $\mu:S^4\to \R^2$ whose origami polytope is the union of two copies of the triangle, $\mu(S^4)=\Delta=\Delta_1\cup\Delta_2$.  
The open covering $\{\U', M_{i,k}^{(j)}\}_{i,j,k}$ consists of the inverse images of the following two copies of 5 open subsets of $\Delta_1$ ($=\Delta_2$) for any small $\vep >0$ : 
\begin{itemize}
\item $\tilde\Delta_Z$ : small open neighbourhood of the hypotenuse $\xi_1+\xi_2=1/2$. 
\item $\tilde\Delta^{(0)}_1=\tilde\Delta^{(0)}_2$ : small open ball of radius $\vep >0$ centered at $(0,0)$. 
\item $\tilde\Delta^{(1)}_{1,1}=\tilde\Delta^{(1)}_{2,1}$ : small open neighbourhood of the line segment,  $0\leq \xi_1<\vep, \vep/2\leq \xi_2 \leq 1-\vep$. 
\item $\tilde\Delta^{(1)}_{1,2}=\tilde\Delta^{(1)}_{2,2}$ : small open neighbourhood of the line segment $0\leq \xi_2<\vep, \vep/2\leq \xi_1 \leq 1-\vep$. 
\item ${\rm int}\Delta_1={\rm int}\Delta_2$. 
\end{itemize}
In this case the cracks consist of the inverse images of two compact subsets ${c_{1,1}^Z}={c_{2,1}^Z}$ and ${c_{1,2}^Z}={c_{2,2}^Z}$ defined by  
$$
{c_{1,1}^Z}\left(={c_{2,1}^Z}\right)  \ : \ \xi_1=0, 1-\vep\leq \xi_2 \leq 1-\vep/2 
$$ and 
$$
{c_{1,2}^Z}\left(={c_{2,2}^Z}\right) \ : \ \xi_2=0, 1-\vep\leq \xi_1\leq 1-\vep/2. 
$$

\begin{figure}[h]
\input{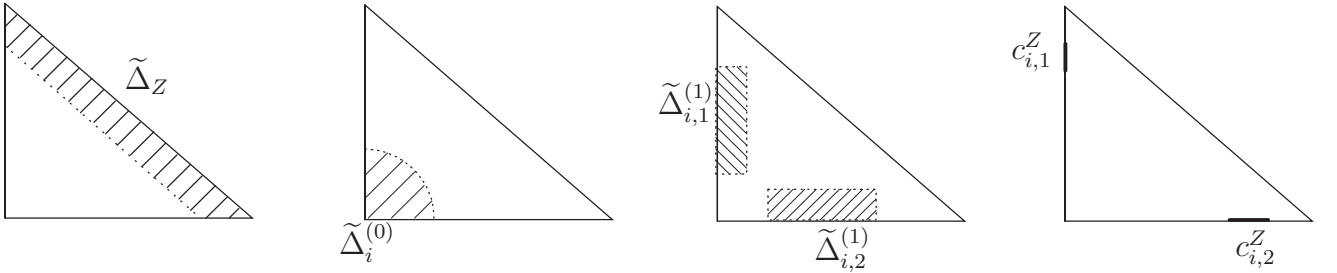}
\caption{Covering of the $S^4$}
\label{Covering of the $S^4$}
\end{figure}

\end{example}

\section{Compatible system on toric origami manifolds}
\label{Compatible system on toric origami manifolds}
In this section we construct a {\it compatible system (of Dirac-type operators)}  on toric origami manifolds. 
The notion of compatible system is introduced in \cite{Fujita-Furuta-Yoshida2}, which is a family of Dirac-type operators along leaves of compatible fibration and satisfies some anti-commutativity. See also Definition~\ref{compatible system}. 

\begin{assump}\label{assump2}
In this section we consider a toric origami manifold $(M,\omega,T,\mu)$ satisfying the following assumption. 
\begin{itemize}
\item $(M,\omega,T,\mu)$ satisfies Assumption~\ref{assump}. 
\item The de Rham cohomology class $[\omega]$ has an integral lift in $H^2(M,\Z)$. 
\item A $T$-equivariant pre-quantizing line bundle $(L,\nabla)$ is fixed. 
Namely, $L$ is a $T$-equivariant Hermitian line bundle over $M$ and $\nabla$ is a $T$-invariant Hermitian connection whose curvature form is equal to $-\sqrt{-1}\omega$. 
\end{itemize}
\end{assump}

Together with the assumptions we may choose a stable almost complex structure $\tilde J$ as in Theorem~\ref{stablealmostcomplex} so that the tangent bundle of each symplectic submanifold $\mu^{-1}({\rm int}(\Delta_{i,k}^{(j)}))$ is preserved by $\tilde J$ for all $i,j$ and $k$. 
Under the above assumption we use the $\Z/2$-graded Clifford module bundle $W_L$ as in the end of Section~\ref{Stable almost complex structure and Clifford module bundle}. 
As it is shown in Section~\ref{Compatible fibration on toric origami manifolds}, $M\setminus \bigcup_{i,k}C_{i,k}^{Z}$ has a structure of good compatible fibration $\{\U', M_{i,k}^{(j)}\}_{i,j,k}$. Since $\{M_{i,k}^{(j)}\}_{i,j,k}$ is a good compatible fibration on an open toric manifold $M\setminus \overline{\U}$, we have a compatible system $\{D_{i,k}^{(j)}\}_{i,j,k}$ on it as in \cite[Theorem~5.1]{Fujita-Furuta-Yoshida2}. Namely for each $i,j,k$ we have the following. 
\begin{itemize}
\item $D_{i,k}^{(j)}$ is a first order formally self-adjoint differential operator of degree-one, which acts on the space of smooth sections of $W_{L}|_{M_{i,k}^{(j)}}$. 
\item $D_{i,k}^{(j)}$ contains only the differentials along the $G_{i,k}^{(j)}$-orbits. 
\item For each $x\in M_{i,k}^{(j)}$, the restriction of $D_{i,k}^{(j)}$ to the orbit $G_{i,k}^{(j)}\cdot x$ is a Dirac-type operator on the $\Z/2$-graded $Cl(T(G_{i,k}^{(j)}\cdot x))$-module bundle $W_L|_{G_{i,k}^{(j)}\cdot x}$. 
\item Let $\tilde u$ be a $G_{i,k}^{(j)}$-invariant section of the normal bundle to the orbit $G_{i,k}^{(j)}\cdot x$. Then $D_{i,k}^{(j)}$ anti-commutes with the Clifford multiplication $c(\tilde u)$ of $\tilde u$ : 
\begin{equation}\label{anti-commutatibity}
D_{i,k}^{(k)}c(\tilde u)+c(\tilde u)D_{i,k}^{(k)}=0.
\end{equation}
\end{itemize}
Now we construct a differential operator $D_{Z}$ along the $S^1$-orbits on $\U$. We first study the product structure of $W|_{\U}$. 
Hereafter we use the identification $\U=Z\times (-\vep, \vep)\cong(Z\times S^1\times(-\vep,\vep))/S^1$ with respect to the diagonal $S^1$-action. 
According to Remark~\ref{stablecomplexstrrem}(2) we may assume that the almost complex structure $\tilde J|_{\U}$ in Theorem~\ref{stablealmostcomplex}  on $T\U\oplus \R^2\cong \pi^*TB\oplus T_{\pi}Z\oplus\R\oplus\R^2$ is the direct sum of almost complex structures on the symplectic vector bundle $\pi^*TB$ and the trivial bundle $T_{\pi}Z\oplus\R\oplus\R^2$ of real rank 4. Then we have 
$$
W|_{\U}={\rm Hom}_{Cl_2}(W_2, \wedge_{\C}^{\bullet}(T\U\oplus\R^2))=\pi^*(\wedge
_{\C}^{\bullet}TB)\otimes{\rm Hom}_{Cl_2}(W_2, \wedge_{\C}^{\bullet}(T_{\pi}Z\oplus\R^3)). 
$$
On the other hand we have the commutative diagram of bundle maps 
\[
\xymatrix{ 
T_{\pi}Z \ar[d] & \ar[l]  {q}^*(T_{\pi}Z)\cong p^*(TS^1)
 \ar[r] \ar[d] & TS^1 \ar[d]\\
Z & \ar[l]^{q}  Z \times S^1 \ar[r]_p & S^1,  } 
\]
where $p:Z\times S^1\to S^1$ is the projection to the $S^1$-factor and $q:Z\times S^1\to (Z\times S^1)/S^1\cong Z$, $(z,t)\mapsto zt^{-1}$ is the quotient map with respect to the diagonal action of $S^1$. 
The isomorphism in the middle column is given by the differential of the map $S^1\to Z$, $t\mapsto zt^{-1}$ for $z\in Z$. The commutative diagram implies that the vector bundle $T_{\pi}Z\oplus\R^3\to\U\cong (Z\times S^1\times(-\vep,\vep))/S^1$ can be obtained as a quotient bundle of $p^*(TS^1)\oplus \R^3\to Z\times S^1\times (-\vep,\vep)$. 
In particular ${\rm Hom}_{Cl_2}(W_2, \wedge_{\C}^{\bullet}(T_{\pi}Z\oplus\R^3))\to\U$ can be obtained as a quotient bundle of ${\rm Hom}_{Cl_2}(W_2, \wedge_{\C}^{\bullet}(TS^1\oplus\R^3))\to Z^1\times S^1\times (-\vep,\vep)$, 
where the complex structure on $TS^1\oplus\R^3$ is given by the same formula for $B_t$ as in the proof of \cite[Theorem~2]{SilvaGuilleminWoodward} under a trivialization. Note that ${\rm Hom}_{Cl_2}(W_2, \wedge_{\C}^{\bullet}(TS^1\oplus\R^3))$ has a structure of $\Z/2$-graded $Cl(TS^1\oplus \R)$-module bundle over $S^1\times(-\vep,\vep)$. 

Now we decompose the line bundle $L$ over $\U$. 
Let $(L_0, \nabla)\to S^1\times(-\vep, \vep)$ be the pre-quantizing line bundle over the folded cylinder as in Appendix~\ref{A computation of local index of the folded cylinder}. 

\begin{prop}\label{product structure of L}
If we take $\vep$ small enough, then 
the diffeomorphism $\varphi:\U\stackrel{\cong}{\to} Z\times (-\vep,\vep)$ as in Theorem~\ref{Morsermodel} can be lifted to an isomorphism between $L|_{\U}\to \U$ and $(L|_Z\boxtimes L_0)/S^1\to (Z\times S^1\times(-\vep,\vep))/S^1=Z\times(-\vep, \vep)$.  
\end{prop}
\begin{proof}
Note that there exists the canonical isomorphism $\tilde\varphi_0$ between $\iota_Z^*L$ and $\iota_0^*\left((L|_Z\boxtimes L_0)/S^1\right)$. Fix a Hermitian connection of $(L|_Z\boxtimes L_0)/S^1$. Then the we have the required isomorphism by using $\tilde\varphi_0$ and the parallel transport. 
\end{proof}

Summarising we have the following. 

\begin{prop}\label{product structure of Cl}
Let $W_{B,L_B}:=\wedge^{\bullet}_{\C}TB\otimes (L|_Z/S^1)$ be a $\Z/2$-graded $Cl(TB)$-module bundle over $B$.  Let $W_{0,L_0}:={\rm Hom}_{Cl_2}(W_2, \wedge_{\C}^{\bullet}(TS^1\oplus\R^3))\otimes L_0$ be a $\Z/2$-graded $Cl(TS^1\oplus\R)$-module bundle over $S^1\times(-\vep,\vep)$ as in the above construction. 
The $\Z/2$-graded Clifford module bundle $W_L|_{\U}\to \U$ is isomorphic to the quotient bundle of the tensor product $\pi^*W_{B,L_B}\otimes p^*W_{0,L_0}\to Z\times S^1\times(-\vep,\vep)$ with respect to the diagonal $S^1$-action, where $\pi:Z\times S^1\times(-\vep,\vep)\to B$ and $p:Z\times S^1\times(-\vep,\vep)\to S^1\times(-\vep,\vep)$ are natural projections. 
\end{prop}

Let $D_{S^1}$ be a Dirac-type operator along the $S^1$-orbits in $S^1\times(-\vep,\vep)$, 
which acts on the space of smooth sections of $W_{0,L_0}$. See Appendix~\ref{A computation of local index of the folded cylinder} for the explicit description of $D_{S^1}$. Let $\epsilon_{B}$ be the map representing the $\Z/2$-grading of $W_{B,L_B}$, i.e., $\epsilon_B(v)=(-1)^{{\rm deg}(v)}(v)$ for $v\in W_{B,L_B}$.  The product of operators $\epsilon_B\otimes D_{S^1}$ is $S^1$-invariant, and it induces a differential operator $D_{Z}$ acting on the smooth sections of $W|_{\U}$ through the isomorphism in Proposition~\ref{product structure of L}. 
Since the $S^1$-action on $Z$ is given by a subgroup of $T$,   $D_{Z}$ is a differential operator along the $S^1$-orbits and satisfies the anti-commutativity as in (\ref{anti-commutatibity}). 

\begin{prop}
The family of differential operators $\{D_{Z}, D_{i,k}^{(j)}\}_{i,j,k}$ is a compatible system on the compatible fibration defined by the torus actions $\{S^1\curvearrowright{\mathcal U}', \ G_{i,k}^{(j)}\curvearrowright M_{i,k}^{(j)}\}_{i,j,k}$. 
\end{prop}

\subsection{Acyclicity of the compatible system}
\label{Acyclicity of the compatible system}
In this section we determine the condition for the compatible system $\{D_Z,D_{i,k}^{(j)}\}_{i,j}$  to be {\it acyclic} (\cite[Definition~6.10]{Fujita-Furuta-Yoshida2} or Definition~\ref{strongly acyclic}).  

Let ${\mathfrak g}_{i,k}^{(j)*}$ be the dual of the Lie algebra of the subtorus $G_{i,k}^{(j)}$ and $({\mathfrak g}_{i,k}^{(j)*})_{\Z}$ the integral weight lattice of ${\mathfrak g}_{i,k}^{(j)*}$. Let $\iota_{i,k}^{(j)}:\g_{i,k}^{(j)}\to \g$ be the inclusion of the Lie subalgebra. Note that the composition $\mu_{i,k}^{(j)}:=(\iota_{i,k}^{(j)*})\circ \mu : M_{i,k}^{(j)} \to {\mathfrak g}_{i,k}^{(j)*}$ is the moment map for the Hamiltonian $G_{i,k}^{(j)}$-action on $M_{i,k}^{(j)}$. We put $M_{i,k}^{(j)\circ}:=M_{i,k}^{(j)}\setminus (\mu_{i,k}^{(j)})^{-1}((\g_{i,k}^{(j)*})_{\Z})$.

\begin{prop}\label{acyclic1} 
For each $x\in M_{i,k}^{(j)\circ}$, we have $\ker(D_{i,k}^{(j)}|_{G_{i,k}^{(j)}\cdot x})=0$.  
\end{prop}
\begin{proof}
Note that for each $x\in M_{i,k}^{(j)}$ the kernel of $D_{i,k}^{(j)}|_{G_{i,k}^{(j)}\cdot x}$ vanishes if and only if there are no non-trivial global parallel sections of $L|_{G_{i,k}^{(j)}\cdot x}$. The proposition follows from the fact that if there exists a global parallel section, then we have $\mu_{i,k}^{(j)}(x)=\iota_{i,k}^{(j)*}(\mu(x))$ lies in the integral weight lattice $({\mathfrak g}_{i,k}^{(j)*})_{\Z}$. 
\end{proof}

We may take $\vep>0$ small enough so that $\mu(\U)=\mu(Z\times (-\vep,\vep))$ does not contain any integral lattice points outside $\mu(Z)=\Delta_Z$. Then we have the following by the same argument as that for Proposition~\ref{acyclic1}. 

\begin{prop}\label{acyclic2} 
For each $x\in \U'\setminus Z$, we have $\ker(D_Z|_{S^1\cdot x})=0$.  
\end{prop}

We put $V:=\left(\U'\cup\bigcup_{i,j,k}M_{i,k}^{(j)\circ}\right)\setminus \left(Z\cup\bigcup_{i,k}C_{i,k}^{Z}\right)$. Then $M\setminus V$ is compact. Since $\{S^1\curvearrowright{\mathcal U}', \ G_{i,k}^{(j)}\curvearrowright M_{i,k}^{(j)}\}_{i,j,k}$ is a good compatible fibration one can see that the following four types of the anti-commutators on the intersections are non-negative. 
\begin{itemize}
\item $D_{i,k}^{(j)}D_{i,k'}^{(j')}+D_{i,k'}^{(j')}D_{i,k}^{(j)}$ on $M_{i,k}^{(j)}\cap M_{i,k'}^{(j')}$, 
and 
\item  $D_ZD_{i,1}^{(n)}+D_{i,1}^{(n)}D_Z$ on $\U'\cap M_{i,1}^{(n)}$. 
\end{itemize}
See \cite[Proposition~5.8, Lemma~5.9]{Fujita-Furuta-Yoshida2} for example. Together with Proposition~\ref{acyclic1} this fact implies the following.  
\begin{prop}\label{acyclic3} 
The compatible system $\{D_{Z}, D_{i,k}^{(j)}\}_{i,j,k}$ is acyclic over $V$. 
\end{prop}

\subsection{Localization formula and Danilov-type formula}
\label{Localization formula and Danilov type formula}
As in Definition~\ref{origamiRR}, the Riemann-Roch number $RR(M,L)$ is defined for any origami manifold $(M,\omega)$ with pre-quantizing line bundle $(L,\nabla)$. If $(M,\omega)$ is a toric origami manifold with the action of a torus $T$, then the resulting index is an element of the character ring $R(T)$ of $T$. In this case we call the index the {\it equivariant Riemann-Roch number} or {\it Riemann-Roch character} and is denoted by $RR_T(M,L)$. 

We use notations in the previous sections and assume Assumption~\ref{assump2}. For each $i,j(\neq n),$ and $k$ we may assume that 
$$
\tilde\Delta_{i,k}^{(j)}\cap{\rm int}\Delta_i\cap \Tt_{\Z}^*=\emptyset,   
$$ and we take and fix a $T$-invariant small open neighbourhood $V_{i,k}^{(j)}$ of $(\mu_{i,k}^{(j)})^{-1}((\g_{i,k}^{(j)*})_{\Z})$ for each $i,j$ and $k$. By the above assumption one has that if $j\neq n$,  then  $V_{i,k}^{(j)}\cap \mu^{-1}(\Tt_{\Z}^*)$ consists of the inverse image of lattice points in the boundary $\partial \Delta_i=\Delta_i\setminus {\rm int}\Delta_i$. 
We also take and fix a small open neighbourhood $V_{i,k}^Z$ of the crack $C_{i,k}^Z$ so that it does not contain any integral points for each $i$ and $k$. Note that each open subset $V_{i,k}^{(j)}\cap V$ (resp. $V_{i,k}^{Z}\cap V$) with compact complement $V_{i,k}^{(j)}\setminus V_{i,k}^{(j)}\cap V=(\mu_{i,k}^{(j)})^{-1}((\g_{i,k}^{(j)*})_{\Z})(\supset M_{i,k}^{(j)}\cap \mu^{-1}(\g_{\Z}^*))$ (resp. $V_{i,k}^{Z}\setminus V_{i,k}^{Z}\cap V=C_{i,k}^{Z})$ is equipped with an acyclic compatible system by Proposition~\ref{acyclic3}, and hence,  the $T$-equivariant local index $\ind_T(V_{i,k}^{(j)}, V_{i,k}^{(j)}\cap V)$ (resp. $\ind_T(V_{i,k}^Z,V_{i,k}^Z\cap V)$) is defined (Theorem~\ref{def of local ind}) . As in the same way one can define the $T$-equivariant local index for the fold, $\ind_T(\U',\U'\setminus Z)$, is defined.

The localization formula (Theorem~\ref{localizationprototype}) implies that the Riemann-Roch character is localized at $\mu^{-1}(\g_{\Z}^*)\cup Z\cup \bigcup_{i,k}C_{i,k}^Z\subset M\setminus V$ as follows. 
\begin{theorem}\label{localization of RR}
Under Assumption~\ref{assump2} we have the localization formula of $T$-equivariant index 
$$
RR_T(M,L)=\ind_T(\U',\U'\setminus Z)+\sum_{i,j,k}\ind_T(V_{i,k}^{(j)}, V_{i,k}^{(j)}\cap V)+\sum_{i,k}\ind_T(V_{i,k}^Z, V_{i,k}^Z\cap V). 
$$
\end{theorem}
By computing the contributions $\ind_T(\U',\U'\setminus Z)$ 
(Theorem~\ref{foldind=0}), $\ind_T(V_{i,k}^{(j)}, V_{i,k}^{(j)}\cap V)$ 
(Theorem~\ref{positive contribution}, Theorem~\ref{negative contribution}) 
and $\ind_T(V_{i,k}^Z, V_{i,k}^Z\cap V)$ (Theorem~\ref{crackind=0}) in the subsequent section, we have the following Danilov-type formula. 

\begin{theorem}\label{origamiDanilov}
Under Assumption~\ref{assump2} we have the following equality as elements in the character ring $R(T)$.\begin{equation}\label{origamiDanilovformula}
RR_T(M,L)=\sum_{\xi_+\in \mu(M^+)\cap\Tt^*_{\Z}}\C_{(\xi_+)}-\sum_{\xi_-\in \mu(M^-)\cap\Tt^*_{\Z}}\C_{(\xi_-)},  
\end{equation}
where for each $\xi \in \Tt_{\Z}^*$ we denote by $\C_{(\xi)}$ the irreducible representation of $T$ whose weight is given by $\xi$.  
\end{theorem}

To compute the local contributions in the subsequent sections, we will use the following notations. 
We divide the collection of Delzant polytopes $\{\Delta_i\}_{i=1,\cdots,N}$ into two subsets, 
$$
\{\Delta_i\}_{i=1,\cdots, N}=\{\Delta_{i}^+\}_{i=1,\cdots, N_+}\cup\{\Delta_{i}^-\}_{i=1,\cdots, N_-}, 
$$where $N_++N_-=N$ and the sign is determined by the condition $\displaystyle\mu(M^{\pm})=\bigcup_{i=1}^{N_{\pm}}\Delta_i^{\pm}$. In a similar way we also use notations $\Delta_{i,k}^{(j)\pm}$,  $V_{i,k}^{(j)\pm}$, $\mu_{i,k}^{(j)\pm}$ and $\g_{i,k}^{(j)\pm}$. 


In terms of this notations the formula (\ref{origamiDanilovformula}) can be rewritten as 
$$
RR_T(M,L)=\sum_{i,j,k}\left(\sum_{\xi_+\in {\rm int}\Delta_{i,k}^{(j)+}\cap\Tt^*_{\Z}}\C_{(\xi_+)}-\sum_{\xi_-\in {\rm int}\Delta_{i,k}^{(j)-}\cap\Tt^*_{\Z}}\C_{(\xi_-)}\right)
$$


\begin{example}\label{exsphere6}
Consider the toric origami manifold $(S^{2n},\omega)$, the unit sphere, with the moment map $\mu:S^{2n}\to \R^n$ as in Example~\ref{exsphere3}, whose origami polytope is the union of two copies of the $n$-simplex, $\mu(S^{2n})=\Delta=\Delta_1\cup\Delta_2$. Since $\mu((S^{2n})^+)\cap\Tt^*_{\Z}=\mu((S^{2n})^-)\cap\Tt^*_{\Z}$,  one has 
$RR_T(S^{2n}, L)=0$ for any $T$-equivariant pre-quantizing line bundle $L$. 

Note that if we use the folded symplectic form $k\omega$ for any positive constant $k$, then the origami polytope for $(S^{2n}, k\omega)$ is the similar extension with ratio $k$ of the original origami polytope. In this case one also has $RR_T(S^{2n},L_k)=0$ for any $T$-equivariant pre-quantizing line bundle $L_k$. 
\end{example}

\subsection{Comments on another possible approaches}
The formula (\ref{origamiDanilovformula}) in Theorem~\ref{origamiDanilov} itself can be obtained as a consequence of the cobordism theorem \cite[Theorem~4.1]{Silva-Guillemin-Pires} and Danilov's theorem for symplectic toric manifolds. 

There is an another possible approach which uses the theory of {\it multi-fans} introduced in \cite{Hattori-Masuda1}. The equivariant index formula \cite[Theorem~11.1]{Hattori-Masuda1}, which is based on the fixed point formula, would be available to the left hand side of (\ref{origamiDanilovformula}). In fact as it is shown in \cite{Masuda-Park} one can associate a multi-fan for each oriented toric origami manifold.  

It would be possible to show the formula (\ref{origamiDanilovformula}) by using the theory of {\it transverse index} in \cite{Braverman}\cite{ParadanVergne}. In \cite{Braverman} it was shown that the Riemann-Roch character $RR_T(M,L)$ can be realized as a perturbation of Dirac operator by the Clliford multiplication of the Kirwan vector field of the moment map. By considering the perturbation $RR_T(M,L)$ is localized at the zero locus of the Kirwan vector field, i.e., the fixed point set $M^{T}$. Under Assumption~\ref{assump}, the fold has a free $S^1$-action, and hence, there are no contributions of the fold to $RR_T(M,L)$. In particular $RR_T(M,L)$ is the sum of contributions of the vertices of the image of the moment map $\mu(M\setminus Z)$.  As in \cite[Example~13]{Vergne3} the contribution from a fixed point is infinite sum of one dimensional representations of $T$ in general. It implies that $RR_T(M,L)$ is expressed as a cancellation of infinite sum of one dimensional representations. See also \cite{HajimeS1}  for the infinite dimensional nature of the transverse index and the finite dimensional nature of the index theory in \cite{Fujita-Furuta-Yoshida1, Fujita-Furuta-Yoshida2}.  

In contrast to these approaches our proof is direct and geometric, which detects the contribution of each lattice point directly and contains a new proof of original Danilov's theorem as a special case.  
\section{Computation of the local contribution} 
\label{Computation of the local contribution}
\subsection{Toric case} 
\label{Toric case}
In this subsection we consider the symplectic toric case, i.e., toric origami manifolds with empty fold. We first summarize the set-up and notations. 

Let $X$ be a $2n$-dimensional symplectic manifold equipped with a Hamiltonian torus action of an $n$-dimensional torus $G$. We assume that there exists a $G$-equivariant pre-quantizing line bundle $L_X\to X$.  Let $\mu_X:X\to {\mathfrak g}^*={\rm Lie}(G)^*$ and $\Delta_X=\mu_X(X)$ be the corresponding moment map and the Delzant polytope.   We take and fix an $m$-dimensional face $\Delta'$ of $\Delta_X$ and a point $\xi$ in the relative interior ${\rm int}(\Delta')$. Let $F:=\mu_X^{-1}(\xi)$ be the $m$-dimensional isotropic torus in $X$ and $X':=\mu_X^{-1}(\Delta')$ be the $2m$-dimensional symplectic submanifold of $X$. We take and fix a point $x\in F\subset X'$. Let $H$ be the stabilizer subgroup at $x$ with respect to $G$-action and $H^{\perp}$ the complementary orthogonal subtorus of $H$ in $G$ with respect to a rational metric of $\g$. Note that $H$ (resp. $H^{\perp}$) is  an $n-m$-dimensional (resp. $m$-dimensional) subtorus of $G$. We denote the inclusion map of Lie-algebra and its dual by $\iota_{H}:{\rm Lie}(H)=\h\to \g$ and  $\iota_{H}^*:\g^*\to\h^*$ respectively.

We first give following comments. 
\begin{itemize}
\item Since the computation is purely local, we do not need the compactness of $\Delta_X$.  In fact we only use a part of the Delzant condition near $\xi$.    
\item We fix a $G$-invariant $\omega$-compatible almost complex structure on $X$ so that it also induces a $G$-invariant $\omega$-compatible almost complex structure on the inverse image of each face of $\Delta_X$. 
\item $F$ is a Lagrangian torus in the symplectic submanifold $X'$. 
\item $F$ can be described as the orbit $F=G\cdot x=H^{\perp}\cdot x$. 
\item The intersection $H\cap H^{\perp}$ is a finite Abelian group. 
\item Since $x$ is a fixed point with respect the $H$-action, the moment map image $(\iota_H^*\circ\mu)(x)=\iota_H^*(\xi)$ of $x$ with respect to the $H$-action is an element in the weight lattice $\h_{\Z}^*$. 
\item The argument below still holds when there exists a finite subgroup of $G$ which acts trivially on $X$. In fact in the proof of Lemma~\ref{locinddisc} we deal with the symplectic toric manifold $X_1$ for which such a subgroup $H\cap H_1\cap H_1^{\perp}$ may exist. 
\end{itemize}

If $Y$ is a smooth manifold and $Y'$ is its smooth submanifold, then we denote the normal bundle of $Y'$ in $Y$ by $\nu_Y(Y')$. We also denote the fiber at $y\in Y'$ by $\nu_Y(Y')_y$.  
There exists a $G$-invariant tubular neighbourhood $N_F$ of $F$ and $G$-equivariant diffeomorphism 
$$
N_F\cong(\nu_X(F)_x\times G)/H
=(\nu_X(F)_x\times H^{\perp})/{H\cap H^{\perp}}, 
$$where we use the $G$-action on the right hand side through the identification 
$G=H\cdot H^{\perp}=(H\times H^{\perp})/H\cap H^{\perp}$ arising from the exact sequence 
\[
\begin{split}
H\cap H^{\perp}\to H\times H^{\perp}\to H\cdot H^{\perp}=G \\ 
h\mapsto (h,h^{-1}), (h_1, h_2)\mapsto h_1h_2. 
\end{split}
\]
Since $F$ is a Lagrangian torus in $X'$ we have 
$$
\nu_X(F)_x\times H^{\perp}=
\nu_X(X')_x\times\nu_{X'}(F)_x\times H^{\perp}
=\nu_X(X')_x\times T^*_x(H^{\perp}\cdot x)\times H^{\perp}
=\nu_X(X')_x\times T^*H^{\perp},
$$ and hence, we have a $G$-equivariant isomorphism 
\begin{equation}\label{nbdofF}
N_F\cong (\nu_X(X')_x\times T^*H^{\perp})/H\cap H^{\perp}. 
\end{equation}

Now we describe the restriction $L_X|_{N_F}$. We first define an $H$-equivariant line bundle $L_1:=\nu_X(X')_x\times L_X|_x\to \nu_X(X')_x$, where we regarded $L_X|_x$ as a representation of $H$. Note that $\nu_X(X')_x$ has a natural symplectic structure and $L_1$ is equipped with a structure of pre-quantizing line bundle with respect to the symplectic structure. 
Let $L_2$ be the pull-back of $L_X|_{N_F}$ with respect to the natural map 
$T^*H^{\perp}\to (\nu_X(X')_x\times T^*H^{\perp})/H\cap H^{\perp}$, which is an $H\times H^{\perp}$-equivariant line bundle over $T^*H^{\perp}$. Note that though $H$-action on $T^*H^{\perp}$ is trivial, the action on $L_2$ is non-trivial in general.  We define an $H^{\perp}$-equivariant line bundle $\hat L_2\to T^*H^{\perp}$ by $\hat L_2:={\rm Hom}(L_X|_x,L_2)$. Then $\hat L_2$ is isomorphic to $L_2$ as $H^{\perp}$-equivariant line bundle and the induced $H$-action on $\hat L_2$ is trivial. We have two line bundles with connection $(L_1\boxtimes \hat L_2)/H\cap H^{\perp}$ and $L_X|_{N_F}$ over $(\nu_X(X')_x\times T^*H^{\perp})/H\cap H^{\perp}=N_F$. The restrictions of these two line bundles to the zero-section $F$ in $N_F$ are isomorphic to each other as line bundles with connection. The Darboux type theorem (\cite[Proposition~7.11]{Fujita-Furuta-Yoshida3}) implies that the $G$-equivariant isomorphism can be extended to a $G$-invariant neighbourhood of $F$.

\begin{remark}
Strictly speaking we have to consider the data on sufficiently small neighbourhoods of the origin in $\nu_X(X')_x$ and the zero section $H^{\perp}$ in $T^*H^{\perp}$ as a  Lagrangian torus to consider the above isomorphisms and  the local indices  in the subsequent argument, though, we use the same notations $\nu_X(X')$ and $T^*H^{\perp}$ to simplify the notations. 
\end{remark}

Let $\Delta_1, \ldots, \Delta_{n-m}$ be codimension one faces of $\Delta_X$ such that $\Delta'$ is the intersection of them,  $\Delta'=\Delta_1\cap \cdots \cap \Delta_{n-m}$. For each $l=1,2,\ldots, n-m$, let $H_l$ be the circle subgroup of $H$ which acts trivially on the symplectic submanifold $X_l:=\mu_X^{-1}(\Delta_l)$ and $H_l^{\perp}$ the orthogonal complement of $H_l$. If we choose any members $\Delta_{l_1}, \ldots, \Delta_{l_\alpha}$, then we have a locally free action of the intersection $H_{l_1}^{\perp}\cap \cdots \cap H_{l_\alpha}^{\perp}=(H_{l_1}\cdot\cdots \cdot	H_{l\alpha})^{\perp}$ on a small neighbourhood of the inverse image of the complement of a neighbourhood of the boundary $\partial(\Delta_{l_1}\cap\cdots\cap\Delta_{l_\alpha})$ in $\Delta_{l_1}\cap\cdots\cap\Delta_{l_\alpha}$. Such a family of torus actions determines a good compatible fibration as in Section~\ref{Good compatible fibration on toric origami manifolds}.  For each $H_{l}$ we have the decomposition $H_l^{\perp}=(H\cap H_l^{\perp})\cdot H^{\perp}$.
On the other hand there exists a natural action of the product $H\times H^{\perp}$ on $N_F$ under the identification (\ref{nbdofF}). 
Then the above good compatible fibration is induced from the action of the subgroup $(H\cap H^{\perp}_l)\times H^{\perp}$ in $H\times H^{\perp}$. 

The $G$-equivariant local index $\ind_G(N_F,N_F\setminus F)$ is defined by using these structures and it is equal to the $H\cap H^{\perp}$-invariant part of the $H\times H^{\perp}$-equivariant local index $\ind_{H\times H^{\perp}}(\nu_X(X')_x\times T^*H^{\perp}, \nu_X(X')_x\times T^*H^{\perp}\setminus \{0\}\times H^{\perp})$. 
%
For simplicity we use the following type of notations for the equivariant local indices: 
$$
{RR}_H(\nu_X(X')_x):=\ind_{H}(\nu_X(X')_x, \nu_X(X')_x\setminus \{0\})
$$and 
$$
RR_{H^{\perp}}(T^*H^{\perp}):=\ind_{H^{\perp}}(T^*H^{\perp}, T^*H^{\perp}\setminus H^{\perp}). 
$$

\begin{lemma}\label{locinddisc}
${RR}_H(\nu_X(X')_x)=\C_{(\iota_H^*(\xi))}=L_X|_x\in R(H)$.
\end{lemma}
\begin{proof}
The second equality follows from the property of the moment map and the Kostant formula. 
We show the first equality by induction on  $n-m=\dim(\nu_X(X')_x)/2$. If $n-m=1$, then the equality follows from the direct computation. See \cite[Example~2.3]{yoshidasymplecticcut} for example. Suppose that $n-m$ is grater than 1 and the statement holds for any situation with codimension $n-m-1$. We consider the decomposition $\nu_X(X')=\nu_{X}(X_1)\oplus\nu_{X_1}(X')$ and $H=H_1\cdot(H\cap H_1^{\perp})$. According to the decomposition the $H$-action on $\nu_X(X')$ factors the action of the product of $H_1$-action on $\nu_X(X_1)$ and $H\cap H_1^{\perp}$-action on $\nu_{X_1}(X')$. By Proposition~\ref{cobandprod} we have that $RR_{H}(\nu_X(X')_x)$ is equal to the $H_1\cap(H\cap H_1^{\perp})$-invariant part of the product $RR_{H_1}(\nu_X(X_1)_x)\otimes RR_{H\cap H_1^{\perp}}(\nu_{X_1}(X')_x)$. Note that $H_1^{\perp}$-action on $X_1$ gives a structure of a symplectic toric manifold whose momentum polytope is $\iota^*_{H^{\perp}_1}(\Delta_1)$. By the assumption of the induction we have $RR_{H_1^{\perp}}(\nu_{X_1}(X')_x)=\C_{(\iota_{H_1^{\perp}}^*(\xi))}$. 
By considering the subgroup $H\cap H_1^{\perp}$ we have $RR_{H\cap H_1^{\perp}}(\nu_{X_1}(X')_x)=\C_{(\iota_{H\cap H_1^{\perp}}^*(\xi))}$, and hence,  
$$
RR_{H_1}(\nu_X(X_1)_x)\otimes RR_{H\cap H_1^{\perp}}(\nu_{X_1}(X')_x)=\C_{(\iota_{H_1}^*(\xi))}\otimes\C_{(\iota_{H\cap H_1^{\perp}}^*(\xi))}=\C_{(\iota_{H_1}^*(\xi)\oplus\iota_{H\cap H_1^{\perp}}^*(\xi))}. 
$$Note that under the natural isomorphism ${\rm Lie}(H_1)^*\oplus {\rm Lie}(H\cap H_1^{\perp})^*  \cong \h^*$ we have $\iota_{H_1}^*(\xi)\oplus\iota_{H\cap H_1^{\perp}}^*(\xi)=\iota_H^*(\xi)$. As we noted in the beginning of this section $\iota_H^*(\xi)$ is an element of the weight lattice $\h_{\Z}^*$, the $H_1\times (H\cap H_1^{\perp})$-representation $\C_{\iota_{H_1}^*(\xi)\oplus\iota_{H\cap H_1^{\perp}}^*(\xi)}$ induces an $H$-representation $\C_{(\iota_H^*(\xi))}$, and hence,  it implies that $RR_{H_1}(\nu_X(X_1)_x)\otimes RR_{H\cap H_1^{\perp}}(\nu_{X_1}(X')_x)$ decsends to an $H$-representation. 
In particular 
the index $RR_{H_1}(\nu_X(X_1)_x)\otimes RR_{H\cap H_1^{\perp}}(\nu_{X_1}(X')_x)$ is $H_1\cap H\cap H_1^{\perp}$-invariant, and we complete the proof. 
\end{proof}

Let $\iota_{H^{\perp}}:\h^{\perp}\to \g$ be the inclusion and $\iota^*_{H^{\perp}}$ its dual. We may assume that the moment map image $(\iota_{H^{\perp}}^*\circ\mu)(x)=\iota_{H^{\perp}}^*(\xi)$ of $x$ with respect to the $H^{\perp}$-action is an element in the weight lattice $(\h^{\perp})_{\Z}^*$. Otherwise the compatible system on $T^*H$ is acyclic, and hence, the local index $RR_{H^{\perp}}(T^*H^{\perp})$ is zero. 

\begin{lemma}\label{locindcylinder}
$RR_{H^{\perp}}(T^*H^{\perp})=\C_{(\iota_{H^{\perp}}^*(\xi))}\in R(H^{\perp})$. 
\end{lemma}
\begin{proof}
Since the $H^{\perp}$-action on $T^*H^{\perp}$ is free, the induced good compatible fibration(system) on $TH^{\perp}$ consists of two open subsets, a small open neighbourhood of the zero-section $H^{\perp}$ and its complement. On the other hand by fixing a decomposition $H^{\perp}=(S^1)^{m}$, we have a product structure of compatible fibration and compatible system, where the $S^1$-equivariant data is determined by the inclusion $\iota_i : S^1\hookrightarrow (S^1)^{m}=H^{\perp}$ to the $i$th fator for $i=1,\cdots, m$. 
By applying Proposition~\ref{cobandprod} the local index $RR_{H^{\perp}}(T^*H^{\perp})$ is equal to the product of $RR_{S^1}(T^*S^1)$ defined the structure induced form $\iota_i$'s. Then the lemma follows from the computation of $RR_{S^1}(T^*S^1)$ (See \cite[Proposition~5.3]{Fujita-Furuta-Yoshida3} for example.). 
\end{proof}

Together with the product formula, Lemma~\ref{locinddisc} and Lemma~\ref{locindcylinder} imply the following. 
\begin{prop}\label{product HHperp}
We have the equality 
$$
RR_{H\times H^{\perp}}(\nu_X(X')_x\times T^*H^{\perp})=\C_{(\iota_H^*(\xi)\oplus\iota_{H^{\perp}}^*(\xi))}\in R(H\times H^{\perp}).  
$$ 
\end{prop}

\begin{theorem}\label{toriclocalcontribution}
$\ind_G(N_F,N_F\setminus F)\neq 0$ if and only if $\xi\in\g_{\Z}^*$, 
and if $\xi\in \g_{\Z}^*$, then we have $\ind_G(N_F,N_F\setminus F)=\C_{(\xi)}$. 
\end{theorem}
\begin{proof}
As we explained, $\ind_G(N_F,N_F\setminus F)$ is equal to the $H\cap H^{\perp}$-invariant part of $RR_{H\times H^{\perp}}(\nu_X(X')_x\times T^*H^{\perp})$ which is represented by a one-dimensional representation of $H\times H^{\perp}$. Suppose that the invariant part is non-zero. Then the one-dimensional representation $RR_{H\times H^{\perp}}(\nu_X(X')_x\times T^*H^{\perp})$ descends to a representation of $G$. Since $\h^*\oplus\h^{{\perp}*}$ is isomorphic to $\g^*$ by $\iota_H^*\oplus\iota_{H^{\perp}}^*$, the invariant part is equal to the point $\iota_H^*(\xi)\oplus\iota_{H^{\perp}}^*(\xi)=\xi\in \g_\Z^*$ by Proposition~\ref{product HHperp}. Conversely if $\xi\in \g_{\Z}^*$, then $RR_{H\times H^{\perp}}(\nu_X(X')_x\times T^*H^{\perp})$ represents a point in $\g_{\Z}^*$, and hence, a representation of $G$. In particular we have $\ind_G(N_F,N_F\setminus F)\neq 0$ as the $H\cap H^{\perp}$-invariant part. 
\end{proof}

\begin{definition}
A $G$-orbit $F$ is called a {\it Bohr-Sommerfeld orbit} ({\it BS-orbit} for short) if there exists a non-trivial global parallel section on the restriction $(L_X,\nabla)|_F$. 
\end{definition}

\begin{prop}
A $G$-orbit $F$ is BS-orbit if and only if $\ind_G(N_F,N_F\setminus F)\neq 0$. 
\end{prop}	
\begin{proof}
We fix the decomposition $H^{\perp}=(S^1)^{m}$ as in the proof of Lemma~\ref{locindcylinder}. The computation in \cite[Remark~6.10]{Fujita-Furuta-Yoshida1} says that $RR_{S^1}(T^*S^1)$ is isomorphic to the space of parallel sectoins $\Gamma^{\rm par}(S^1, \iota_i^*\hat L_2|_{S^1})$ for each $i=1,\cdots,m$. By the product structure of $(TH^{\perp}, \hat L_2)$ near the zero section, the space of parallel sections $\Gamma^{\rm par}(H^{\perp}, \hat L_2|_{H^{\perp}})$ is generated by a constant section and isomorphic to the product of $\Gamma^{\rm par}(S^1, \iota_i^*\hat L_2|_{S^1})\cong RR_{S^1}(T^*S^1)$. It implies that $\Gamma^{\rm par}(H^{\perp}, \hat L_2|_{H^{\perp}})$ is isomorphic to $RR_{H^{\perp}}(T^*H^{\perp})$ as $H^{\perp}$-representation. 
Similarly by considering the restriction to the origin, we have that the one-dimensional representation $RR_H(\nu_X(X')_x)$ is isomorphic to $L_X|_x$ as $H$-representation. Then we have that $RR_H(\nu_X(X')_x)\otimes RR_{H^{\perp}}(T^*H^{\perp})$ is isomorphic to $\Gamma^{\rm par}(\{0\}\times H^{\perp}, L_X|_x\otimes \hat L_2|_{H^{\perp}})$ as $H\times H^{\perp}$-representation. 

If $F$ is a BS-orbit, then there exists a non-trivial global parallel section $s_F:F\to L|_F$. By considering the pull-back we have a non-trivial global parallel section $\tilde s_F:H^{\perp}\to L_X|_x\otimes \hat L_2|_{H^{\perp}}$, which is $H\cap H^{\perp}$-invariant, and hence,  
it implies $\ind_G(N_F,N_F\setminus F)\neq 0$. 

Conversely suppose that $\ind_G(N_F,N_F\setminus F)\neq 0$. Then by the isomorphism $RR_H(\nu_X(X')_x)\otimes RR_{H^{\perp}}(T^*H^{\perp}) \cong\Gamma^{\rm par}(\{0\}\times H^{\perp}, L_X|_x\otimes \hat L_2|_{H^{\perp}})$ there exists an $H\cap H^{\perp}$-invariant non-trivial global parallel sections of $L_X|_x\otimes\hat L_2|_{H^{\perp}}$. It induces a non-trivial global parallel section $s_F$ of $L_X|_F$ by the natural map $H^{\perp}\to H^{\perp}\cdot x=F$. 
\end{proof}

When we consider the situation in Section~\ref{Localization formula and Danilov type formula} we have 
$$
\sum_{\xi\in\Delta_i}\ind_T(N_F, N_F\setminus F)=\sum_{j,k}\ind_T(V_{i,k}^{(j)+}, V_{i,k}^{(j)+}\cap V). 
$$
As a particular case we have a proof of Danilov's theorem for symplectic toric manifolds. 
\begin{theorem}\label{Danilov}
If $X$ is a closed symplectic toric manifold with pre-quantizing line bundle $L$, then we have the following equality of the $G$-equivariant Riemann-Roch number. 
$$
RR_G(X,L)=\sum_{\xi\in (\Delta_{X})_{\Z}}\C_{(\xi)}. 
$$
\end{theorem}

\subsection{Contribution from the fold} 
In the subsequent subsections we consider the toric origami case as in Theorem~\ref{origamiDanilov}. In this subsection we compute the contribution from the folded part, $\ind_T(\U',\U'\setminus Z)$. 
\begin{theorem}\label{foldind=0}
We have 
$$
\ind_T(\U',\U'\setminus Z)=0
$$as a $T$-equivariant index. 
\end{theorem}
\begin{proof}
As it is showed in Proposition~\ref{product structure of Cl} and by definition of $D_{Z}$, the acyclic compatible system on $\U'\setminus Z$ has a natural product structure between them on $B$ and $S^1\times(-\vep/2,\vep/2)$, and hence, its local index $\ind_T(\U',\U'\setminus Z)$ is equal to the product of them in the sense of the product formula \cite[Theorem~8.8]{Fujita-Furuta-Yoshida2}. On the other hand the compatible system on $S^1\times (-\vep/2,\vep/2)$ is the one associated with the natural folded structure on it, and it will be shown in Appendix~\ref{A computation of local index of the folded cylinder} that its local index is equal to $0$. See Proposition~\ref{vanishingoffoldedcylinder}. These facts imply $\ind_T(\U',\U'\setminus Z)=0$. 
\end{proof}

\subsection{Contribution from the positive unfolded part}
We compute the contribution from the unfolded part of the positive orientation, $\ind_T(V_{i,k}^{(j)+}, V_{i,k}^{(j)+}\setminus (\mu_{i,k}^{(j)+})^{-1}((\g_{i,k}^{(j)*})_{\Z}))$. Since $V_{i,k}^{(j)+}$ is away from the fold $Z$, the local situation is same as that for the genuine toric case, and hence,  we can apply Theorem~\ref{toriclocalcontribution}. 

\begin{theorem}\label{positive contribution}
We may choose $\tilde\Delta_{i,k}^{(j)+}$ small enough so that $\tilde\Delta_{i,k}^{(j)+}\cap\Tt_{\Z}^*={\rm int}\Delta_{i,k}^{(j)+}\cap\Tt_{\Z}^*$. 
Then we have 
$$
\ind_T(V_{i,k}^{(j)+}, V_{i,k}^{(j)+}\setminus (\mu_{i,k}^{(j)+})^{-1}((\g_{i,k}^{(j)*})_{\Z}))=\sum_{\xi\in{\rm int}\Delta_{i,k}^{(j)+}\cap\Tt_{\Z}^*}\C_{(\xi)}. 
$$
\end{theorem}
\begin{proof}
Since the compatible system $\{D_{Z}, D_{i,k}^{(j)}\}_{i,j,k}$ is acyclic on $V$, the complement of the inverse images of lattice points, the excision formula implies that the $T$-equivariant local index $\ind_T(V_{i,k}^{(j)+}, V_{i,k}^{(j)+}\setminus (\mu_{i,k}^{(j)+})^{-1}((\g_{i,k}^{(j)*})_{\Z}))$ is equal to the sum of contributions of the inverse image of the lattice point which is contained in $V_{i,k}^{(j)+}$. Each inverse image has a neighborhood of the form $N_F$ as in Subsection~\ref{Toric case}, and hence, the contribution of the lattice point $\xi$ is the representation corresponding to the lattice point $\C_{(\xi)}$. 
\end{proof}

\subsection{Contribution from the negative unfolded part}
We compute the contribution from the unfolded part of the negative orientation, $\ind_T(V_{i,k}^{(j)-}, V_{i,k}^{(j)-}\setminus (\mu_{i,k}^{(j)-})^{-1}((\g_{i,k}^{(j)*})_{\Z}))$. The situation is same as that for the positive unfolded part up to the orientation. 
The difference appears only in the $\Z/2$-grading of the Clifford module bundle. Namely the $\Z/2$-grading in the negative case is opposite to the positive case, and hence, the resulting index has the opposite sign. The proof of the following theorem can be shown by the similar way for the proof of Theorem~\ref{positive contribution}. 

\begin{theorem}\label{negative contribution}
We may choose $\tilde\Delta_{i,k}^{(j)-}$ small enough so that $\tilde\Delta_{i,k}^{(j)-}\cap\Tt_{\Z}^*={\rm int}\Delta_{i,k}^{(j)-}\cap\Tt_{\Z}^*$. 
Then we have 
$$
\ind_T(V_{i,k}^{(j)-}, V_{i,k}^{(j)-}\setminus (\mu_{i,k}^{(j)-})^{-1}((\g_{i,k}^{(j)*})_{\Z}))=-\sum_{\xi\in{\rm int}\Delta_{i,k}^{(j)-}\cap\Tt_{\Z}^*}\C_{(\xi)}. 
$$
\end{theorem}

\subsection{Contribution from the crack}\label{Contribution from the crack}
We compute the contribution from the crack, $\ind_T(V_{i,k}^{Z}, V_{i,k}^{Z}\cap V)=\ind_T(V_{i,k}^{Z}, V_{i,k}^{Z}\setminus C_{i,k}^Z)$, and show that it is equal to 0. Note that each $V_{i,k}^Z$ has two components $V_{i,k}^Z\cap M^+$ and $V_{i,k}^Z\cap M^-$.  Then the open subsets $V_{i,k}^Z\cap M^+$ and $V_{i,k}^Z\cap M^-$ are isomorphic to each other as open symplectic toric manifolds up to their orientations. 
\begin{theorem}\label{crackind=0} 
We have the equality 
$$
\ind_T(V_{i,k}^{Z}, V_{i,k}^{Z}\setminus C_{i,k}^Z)=0
$$ as $T$-equivariant indices for each $i$ and $k$.  
\end{theorem}
\begin{proof}
Since $V_{i,k}^Z\cap M^+$ and $V_{i,k}^Z\cap M^-$ are isomorphic up to their orientations we have 
$$
\ind_T(V_{i,k}^{Z}, V_{i,k}^{Z}\setminus C_{i,k}^Z)=\ind_T(V_{i,k}^{Z}\cap M^+, V_{i,k}^{Z}\cap M^+\setminus C_{i,k}^{Z})+\ind_T(V_{i',k'}^{Z}\cap M^-, V_{i,k}^{Z}\cap M^-\setminus C_{i,k}^{Z})=0. 
$$
\end{proof}

\appendix

\section{Acyclic compatible systems and their local indices}
\label{Acyclic compatible systems and their local indices} 
In this appendix we give a brief summary of some definitions 
of compatible fibration, acyclic compatible system and 
their local indices following \cite{Fujita-Furuta-Yoshida2, Fujita-Furuta-Yoshida3} and \cite{Fujitacobinv}. 
We adopt combinations of definitions in \cite{Fujita-Furuta-Yoshida2} and \cite{Fujita-Furuta-Yoshida3}. 
Let $V$ be a smooth  manifold. 
\begin{definition}\label{compatible fibration}
A {\it compatible fibration on $V$} is a collection of the data 
$\{V_{\alpha}, {\SF}_{\alpha}\}_{\alpha\in A}$ consisting of 
an open covering $\{V_{\alpha}\}_{\alpha\in A}$ of $V$ and 
a foliation $\SF_{\alpha}$ on $V_{\alpha}$ with compact leaves 
which  satisfies the following properties.
\begin{enumerate}
\item The holonomy group of each leaf of $\SF_\alpha$ is finite. 
\item\label{correspondence between foliation and pi_alpha}For each $\alpha$ and $\beta$, if a leaf $L\in \SF_\alpha$ has non-empty intersection $L\cap V_\beta\neq \emptyset$, then, $L\subset V_\beta$. 
\end{enumerate}
\end{definition}


\begin{definition}\label{goodcompatifib}
A compatible fibration $\{V_{\alpha}, {\SF}_{\alpha}\}_{\alpha\in A}$ on $V$ is called {\it good} if for all $\alpha$ and $\beta$ with $V_{\alpha}\cap V_{\beta}\neq \emptyset$ the following condition (i) or (ii) holds. 
\begin{itemize}
\item[(i)] For each leaf $L_{\alpha}\in {\SF}_{\alpha}$, there exists a leaf $L_{\beta}\in{\SF}_{\beta}$ such that $L_{\alpha}\subset L_{\beta}$. 
\item[(ii)] For each leaf $L_{\beta}\in {\SF}_{\beta}$, there exists a leaf $L_{\alpha}\in{\SF}_{\alpha}$ such that $L_{\beta}\subset L_{\alpha}$. 
\end{itemize}
\end{definition}

Let $(V,g)$ be a Riemannian manifold, $W$ a $Cl(TV)$-module bundle over $V$. Suppose that $V$ is equipped with a compatible fibration $\{V_{\alpha}, {\SF}_{\alpha}\}_{\alpha\in A}$. 
We impose the following conditions on the Riemannian metric $g$. 
\begin{assumption}\label{assumption for Riemannian metric}
Let $\nu_\alpha=\{ u\in TV_\alpha \mid g(u,v)=0\ \text{for all }v\in T\SF_\alpha \}$ be the normal bundle of $\SF_\alpha$. Then, $g|_{\nu_\alpha}$ is invariant under holonomy, and gives a transverse invariant metric on $\nu_\alpha$.
\end{assumption}

\begin{definition}\label{compatible system}
A {\it compatible system} on $(\{V_{\alpha}, \SF_{\alpha}\}, W)$ is a data $\{ D_{\alpha}\}_{\alpha \in A}$ satisfying the following properties. 
\begin{enumerate}
\item $D_{\alpha}\colon \Gamma (W|_{V_{\alpha}})\to \Gamma (W|_{V_{\alpha}})$ is an order-one formally self-adjoint differential operator.
\item $D_{\alpha}$ contains only the derivatives along leaves of $\SF_{\alpha}$. 
\item $D_{\alpha}$ is a Dirac-type operator along leaves. 
Namely 
the principal symbol of $D_{\alpha}$ is given by the composition of 
the dual of the natural inclusion $\iota_{\alpha}\colon T\SF_{\alpha}\to TV_{\alpha}$ and the Clifford multiplication 
$c\colon T^*\SF_{\alpha}\cong T\SF_{\alpha}\subset TV_{\alpha} \to \End (W|_{V_{\alpha}})$ . 
\item For a leaf $L\in \SF_\alpha$ let $\tilde u\in \Gamma (\nu_\alpha|_L)$ be a section of $\nu_\alpha|_L$ parallel along $L$. 
$\tilde{u}$ acts on $W|_L$ by the Clifford multiplication $c(\tilde{u})$. Then $D_{\alpha}$ and $c(\tilde{u})$ anti-commute each other, i.e. 
\[
0=\{ D_{\alpha},c(\tilde{u}) \}:=
D_{\alpha}\circ c(\tilde{u})+c(\tilde{u})\circ D_{\alpha}
\]
as an operator on $W|_L$. 
\end{enumerate}
\end{definition}

As in \cite[Lemma~3.4]{Fujita-Furuta-Yoshida2} for each leaf $L\in\SF_{\alpha}$ 
we have a small 
open tubular neighbourhood $V_L$ of $L$ and the 
finite covering  $q_L:\tilde V_L\to V_L$ such that the induced foliation on 
$\tilde V_L$ is a bundle foliation with the projection $\pi_L:\tilde V_L\to \tilde U_L$.

\begin{definition}\label{strongly acyclic}
A compatible system $\{ D_\alpha\}_{\alpha \in A}$ on $(\{V_{\alpha}, \SF_{\alpha}\}, W)$ is said to be {\it acyclic} if 
it satisfies the following conditions.  
\begin{enumerate}
\item The Dirac-type operator 
$q^*_{L}D_{\alpha}|_{\pi_L^{-1}(\tilde b)}$ has zero kernel 
for each $\alpha\in A$, leaf $L\in \SF_{\alpha}$ and 
$\tilde b\in \tilde U_{L}$. 
\item If $V_\alpha\cap V_\beta\neq \emptyset$, then the anti-commutator 
$\{D_{\alpha},D_{\beta}\}:=D_{\alpha}D_{\beta}+D_{\beta}D_{\alpha}$ is a non-negative operator on $V_\alpha\cap V_\beta$. 
\end{enumerate}
\end{definition}

As in \cite[Section~5]{Fujita-Furuta-Yoshida2} we can construct such structures, good compatible fibration and compatible system, on Hamiltonian torus manifolds. Though the good compatible fibrations form a nice class, we have to generalize it to treat the product of such structures. 


\begin{definition}\label{tangential}
Suppose that a compact Lie group $G$ 
acts on a Riemannian manifold $V$ in an isometric way. 
Let $\{V_{\alpha}, \SF_{\alpha}\}_{\alpha\in A}$ be a compatible fibration on $V$. 
If the following conditions are satisfied, then we call the compatible fibration a {\it $G$-tangential compatible fibration} (or {\it tangential compatible fibration} for short). 
\begin{itemize}
\item $\{V_{\alpha}\}_{\alpha\in A}$ is a $G$-invariant open covering of $V$. 
\item Each leaf $L$ of $\SF_{\alpha}$ has 
positive dimension for all $\alpha\in A$. 
\item For each leaf $L$ of $\SF_{\alpha}$ there exists some $x\in V_{\alpha}$ 
such that $L$ is contained in the $G$-orbit $G\cdot x$. 
\end{itemize}

A compatible system on a $G$-tangential compatible fibration is 
called {\it $G$-tangential compatible system} 
(or {\it tangential compatible system} for short). 
\end{definition}

Any non-trivial torus action induces a good compatible fibration, which  
is a tangential compatible fibration. 
Moreover the product of two such good compatible fibrations 
is a tangential compatible fibration which is not good in general. 



\begin{theorem}[Theorem~7.2 and Proposition~7.3 in \cite{Fujita-Furuta-Yoshida2},Theorem~3.7 in \cite{Fujitacobinv}]\label{def of local ind}
Suppose that $V$ is an open subset of $M$ whose complement is compact.  
If $V$ is equipped with a $G$-tangential acyclic compatible system $\{V_{\alpha}, \SF_{\alpha}, D_{\alpha}\}_{\alpha\in A}$, then  
we can define the 
local index $$\ind(M, \{V_{\alpha}, \SF_{\alpha}, D_{\alpha}\}_{\alpha\in A},W)=
\ind(M,V,W)=\ind(M,V)\in\Z,$$ which satisfies the excision formula, sum formula and product formula. 
\end{theorem}

Let us briefly recall the definition of the local index $\ind(M,V,W)$.  
Let $D:\Gamma(W)\to \Gamma(W)$ be a Dirac-type operator. 
We consider 
the perturbation $D_t:=D+t\sum_{\alpha\in A}\rho_{\alpha}D_{\alpha}\rho_{\alpha}$ for $t\gg 1$, where $\{\rho_{\alpha}\}_{\alpha\in A}$ is a family of smooth cut-off functions 
which is constant along leaves of $\SF_{\alpha}$ and satisfies 
some estimates as in \cite[Subection~4.1]{Fujita-Furuta-Yoshida2}. 
Such a perturbation $D_t$ gives a Fredholm operator on the space of $L^2$-sections of $W$. The local index $\ind(M,V)$ is defined as the analytic index of $D_t$ for $t\gg 1$. The excision formula implies the following localization formula of Dirac-type operator. 

\begin{theorem}\label{localizationprototype}
Suppose that $M$ is compact without boundary and an open subset $V$ of $M$ is equipped with a $G$-tangential acyclic compatible system. Then the index of any Dirac-type operator $\ind(W)$ is localized at the complement $M\setminus V$. Namely we have 
$$
\ind(W)=\ind(M,V). 
$$
\end{theorem}

\section{A formula of local indices of vector spaces }
\label{A formula of local indeices of vector spaces}

In this appendix we give a formula of equivariant local indices of vector spaces. For $l=1,2$ let $G_l$ be an $m_l$-dimensional torus and $R_l$ an $m_l$-dimensional Hermitian vector space on which $G_l$ acts unitary and effective way. 
We put the following assumptions for $l=1,2$. 

\begin{assump}
$(1)$ A $G_l$-tangential equivariant compatible fibration (Definition~\ref{tangential}) on $R_{l}^{\times}:=R_l\setminus\{0\}$ is given. 
\\
$(2)$ For the compatible fibration in (1), a $G_l$-tangential equivariant acyclic compatible system on $R_l^{\times}$ is given. 
\end{assump}
By the assumption we have two equivariant local indices 
$\ind_{G_1}(R_1,R_1^{\times})$  and $\ind_{G_2}(R_2,R_2^{\times})$. 
Now we fix $\vep>0$ small enough and define two compatible fibrations and acyclic compatible systems on the product $R:=R_1\times R_2$. 

Define two subsets $R'$ and $R''$ of $R$ by 
$$
R':=\{(v_1, v_2)\in R \ | \ |v_1|>\vep, \ |v_2|<\vep\}, 
$$and 
$$
R'':=\{(v_1, v_2)\in R \  | \  |v_1|<\vep, \  |v_2|>\vep \}. 
$$
We consider a structure of $G_1$-tangential (resp. $G_2$-tangential) compatible fibration on $R'$ (resp. $R''$) induced from the first (resp. second) factor. We also define a subset $R_{\infty}$ of $R$ by 
$$
R_{\infty}:=\{(v_1, v_2)\in R \ | \ |v_1|>\vep/2, \ |v_2|>\vep/2\},  
$$which is also equipped with a compatible fibration and compatible system arising from the product structure. 
Then the union $\tilde{R}_{\infty}:=R'\cup R'' \cup R_{\infty}$ gives an open covering of the complement of a compact neighbourhood of the origin of $R$. Note that the above compatible fibration and compatible system define a $G_1\times G_2$-tangential equivariant compatible fibration and acyclic compatible system on $\tilde R_{\infty}$, and hence, we have the equivariant local index 
$$
\ind_{G_1\times G_2}(R, \tilde R_{\infty}). 
$$

\begin{figure}[h]
{\unitlength 0.1in
\begin{picture}( 27.6000, 22.2000)( 23.1000,-26.5000)
%
{\color[named]{Black}{%
\special{pn 8}%
\special{pa 2600 2600}%
\special{pa 2600 600}%
\special{fp}%
\special{sh 1}%
\special{pa 2600 600}%
\special{pa 2580 668}%
\special{pa 2600 654}%
\special{pa 2620 668}%
\special{pa 2600 600}%
\special{fp}%
}}%
%
{\color[named]{Black}{%
\special{pn 8}%
\special{pa 2600 2600}%
\special{pa 5000 2600}%
\special{fp}%
\special{sh 1}%
\special{pa 5000 2600}%
\special{pa 4934 2580}%
\special{pa 4948 2600}%
\special{pa 4934 2620}%
\special{pa 5000 2600}%
\special{fp}%
}}%
\put(50.7000,-26.5000){\makebox(0,0)[lb]{$R_1$}}%
\put(24.6000,-5.6000){\makebox(0,0)[lb]{$R_2$}}%
\put(30.4000,-27.8000){\makebox(0,0)[lb]{$\vep/2$}}%
\put(36.7000,-27.6000){\makebox(0,0)[lb]{$\vep$}}%
\put(23.8000,-17.1000){\makebox(0,0)[lb]{$\vep$}}%
\put(23.1000,-21.9000){\makebox(0,0)[lb]{$\vep/2$}}%
%
{\color[named]{Black}{%
\special{pn 8}%
\special{pa 3150 2590}%
\special{pa 3150 610}%
\special{dt 0.045}%
}}%
%
{\color[named]{Black}{%
\special{pn 8}%
\special{pa 3710 2600}%
\special{pa 3710 620}%
\special{dt 0.045}%
}}%
%
{\color[named]{Black}{%
\special{pn 8}%
\special{pa 2590 2150}%
\special{pa 4880 2150}%
\special{dt 0.045}%
}}%
%
{\color[named]{Black}{%
\special{pn 8}%
\special{pa 2610 1640}%
\special{pa 4900 1640}%
\special{dt 0.045}%
}}%
%
{\color[named]{Black}{%
\special{pn 4}%
\special{pa 4980 2100}%
\special{pa 4480 2600}%
\special{fp}%
\special{pa 4810 2150}%
\special{pa 4360 2600}%
\special{fp}%
\special{pa 4690 2150}%
\special{pa 4240 2600}%
\special{fp}%
\special{pa 4570 2150}%
\special{pa 4120 2600}%
\special{fp}%
\special{pa 4450 2150}%
\special{pa 4000 2600}%
\special{fp}%
\special{pa 4330 2150}%
\special{pa 3880 2600}%
\special{fp}%
\special{pa 4210 2150}%
\special{pa 3760 2600}%
\special{fp}%
\special{pa 4090 2150}%
\special{pa 3710 2530}%
\special{fp}%
\special{pa 3970 2150}%
\special{pa 3710 2410}%
\special{fp}%
\special{pa 3850 2150}%
\special{pa 3710 2290}%
\special{fp}%
\special{pa 4980 2220}%
\special{pa 4600 2600}%
\special{fp}%
\special{pa 4980 2340}%
\special{pa 4720 2600}%
\special{fp}%
\special{pa 4980 2460}%
\special{pa 4840 2600}%
\special{fp}%
\special{pa 4980 1980}%
\special{pa 4810 2150}%
\special{fp}%
\special{pa 4980 1860}%
\special{pa 4690 2150}%
\special{fp}%
\special{pa 4980 1740}%
\special{pa 4570 2150}%
\special{fp}%
\special{pa 4960 1640}%
\special{pa 4450 2150}%
\special{fp}%
\special{pa 4840 1640}%
\special{pa 4330 2150}%
\special{fp}%
\special{pa 4720 1640}%
\special{pa 4210 2150}%
\special{fp}%
\special{pa 4600 1640}%
\special{pa 4090 2150}%
\special{fp}%
\special{pa 4480 1640}%
\special{pa 3970 2150}%
\special{fp}%
\special{pa 4360 1640}%
\special{pa 3850 2150}%
\special{fp}%
\special{pa 4240 1640}%
\special{pa 3730 2150}%
\special{fp}%
\special{pa 4120 1640}%
\special{pa 3710 2050}%
\special{fp}%
\special{pa 4000 1640}%
\special{pa 3710 1930}%
\special{fp}%
\special{pa 3880 1640}%
\special{pa 3710 1810}%
\special{fp}%
\special{pa 3760 1640}%
\special{pa 3710 1690}%
\special{fp}%
}}%
%
{\color[named]{Black}{%
\special{pn 4}%
\special{pa 3150 1230}%
\special{pa 2600 680}%
\special{fp}%
\special{pa 3150 1110}%
\special{pa 2610 570}%
\special{fp}%
\special{pa 3150 990}%
\special{pa 2720 560}%
\special{fp}%
\special{pa 3150 870}%
\special{pa 2840 560}%
\special{fp}%
\special{pa 3150 750}%
\special{pa 2960 560}%
\special{fp}%
\special{pa 3150 630}%
\special{pa 3080 560}%
\special{fp}%
\special{pa 3710 1070}%
\special{pa 3200 560}%
\special{fp}%
\special{pa 3710 1190}%
\special{pa 3150 630}%
\special{fp}%
\special{pa 3710 1310}%
\special{pa 3150 750}%
\special{fp}%
\special{pa 3710 1430}%
\special{pa 3150 870}%
\special{fp}%
\special{pa 3710 1550}%
\special{pa 3150 990}%
\special{fp}%
\special{pa 3680 1640}%
\special{pa 3150 1110}%
\special{fp}%
\special{pa 3560 1640}%
\special{pa 3150 1230}%
\special{fp}%
\special{pa 3440 1640}%
\special{pa 3150 1350}%
\special{fp}%
\special{pa 3320 1640}%
\special{pa 3150 1470}%
\special{fp}%
\special{pa 3200 1640}%
\special{pa 3150 1590}%
\special{fp}%
\special{pa 3710 950}%
\special{pa 3320 560}%
\special{fp}%
\special{pa 3710 830}%
\special{pa 3440 560}%
\special{fp}%
\special{pa 3710 710}%
\special{pa 3560 560}%
\special{fp}%
\special{pa 3710 590}%
\special{pa 3680 560}%
\special{fp}%
\special{pa 3150 1350}%
\special{pa 2600 800}%
\special{fp}%
\special{pa 3150 1470}%
\special{pa 2600 920}%
\special{fp}%
\special{pa 3150 1590}%
\special{pa 2600 1040}%
\special{fp}%
\special{pa 3080 1640}%
\special{pa 2600 1160}%
\special{fp}%
\special{pa 2960 1640}%
\special{pa 2600 1280}%
\special{fp}%
\special{pa 2840 1640}%
\special{pa 2600 1400}%
\special{fp}%
\special{pa 2720 1640}%
\special{pa 2600 1520}%
\special{fp}%
}}%
%
{\color[named]{Black}{%
\special{pn 4}%
\special{pa 5010 1520}%
\special{pa 3140 1520}%
\special{fp}%
\special{pa 5010 1610}%
\special{pa 3140 1610}%
\special{fp}%
\special{pa 5010 1700}%
\special{pa 3140 1700}%
\special{fp}%
\special{pa 5010 1788}%
\special{pa 3140 1788}%
\special{fp}%
\special{pa 5010 1876}%
\special{pa 3140 1876}%
\special{fp}%
\special{pa 5010 1966}%
\special{pa 3140 1966}%
\special{fp}%
\special{pa 5010 2054}%
\special{pa 3140 2054}%
\special{fp}%
\special{pa 5010 2144}%
\special{pa 3140 2144}%
\special{fp}%
\special{pa 5010 1432}%
\special{pa 3140 1432}%
\special{fp}%
\special{pa 5010 1342}%
\special{pa 3140 1342}%
\special{fp}%
\special{pa 5010 1254}%
\special{pa 3140 1254}%
\special{fp}%
\special{pa 5010 1166}%
\special{pa 3140 1166}%
\special{fp}%
\special{pa 5010 1076}%
\special{pa 3140 1076}%
\special{fp}%
\special{pa 5010 988}%
\special{pa 3140 988}%
\special{fp}%
\special{pa 5010 898}%
\special{pa 3140 898}%
\special{fp}%
\special{pa 5010 810}%
\special{pa 3140 810}%
\special{fp}%
\special{pa 5010 720}%
\special{pa 3140 720}%
\special{fp}%
\special{pa 5010 632}%
\special{pa 3140 632}%
\special{fp}%
\special{pa 5010 544}%
\special{pa 3140 544}%
\special{fp}%
}}%
\put(40.3000,-25.1000){\makebox(0,0)[lb]{$R'$}}%
\put(27.6000,-12.3000){\makebox(0,0)[lb]{$R''$}}%
\put(40.3000,-12.4000){\makebox(0,0)[lb]{$R_{\infty}$}}%
\end{picture}}%
\caption{Open covering $\tilde R_{\infty}$.}
\label{ldomein}
\end{figure}
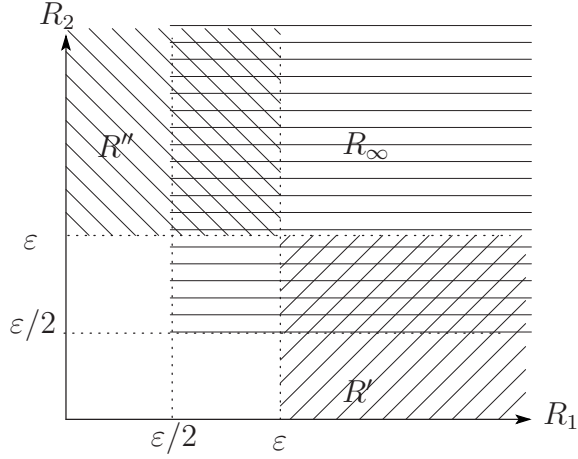

For $l=1,2$ define open  subsets $R_{l,0}$ and $R_{l,\infty}$ of $R_l$ by 
$$
R_{l,0}:=\{v\in R_l \ | \ |v|<\vep\}, 
$$and 
$$
R_{l,\infty}:=\{v\in R_l \ | \ |v|>\vep/2|\}. 
$$We set $R_{\infty}^{\rm prod}:=(R_{1,\infty}\times R_{2,0})\cup (R_{1,0}\times R_{2,\infty})\cup (R_{1,\infty}\times R_{2,\infty})$, which gives an open covering of a complement of a compact neighbourhood of the origin of $R$.
We consider the trivial fibration on $R_{l,0}$ and the $G_l$-tangential compatible fibration on $R_{l,\infty}$. The product of these structures induces a $G_1\times G_2$-tangential equivariant compatible fibration and acyclic compatible system on $R_{\infty}^{\rm prod}$, and hence, we have the equivariant local index 
$$
\ind_{G_1\times G_2}(R, R_{\infty}^{\rm prod}). 
$$

\begin{figure}[h]
{\unitlength 0.1in
\begin{picture}( 27.2000, 22.5000)( 23.1000,-26.7000)
%
{\color[named]{Black}{%
\special{pn 8}%
\special{pa 2600 2600}%
\special{pa 2600 600}%
\special{fp}%
\special{sh 1}%
\special{pa 2600 600}%
\special{pa 2580 668}%
\special{pa 2600 654}%
\special{pa 2620 668}%
\special{pa 2600 600}%
\special{fp}%
}}%
%
{\color[named]{Black}{%
\special{pn 8}%
\special{pa 2600 2600}%
\special{pa 5000 2600}%
\special{fp}%
\special{sh 1}%
\special{pa 5000 2600}%
\special{pa 4934 2580}%
\special{pa 4948 2600}%
\special{pa 4934 2620}%
\special{pa 5000 2600}%
\special{fp}%
}}%
\put(50.3000,-26.7000){\makebox(0,0)[lb]{$R_1$}}%
\put(24.7000,-5.5000){\makebox(0,0)[lb]{$R_2$}}%
\put(30.4000,-28.0000){\makebox(0,0)[lb]{$\vep/2$}}%
\put(36.7000,-27.6000){\makebox(0,0)[lb]{$\vep$}}%
\put(24.1000,-17.1000){\makebox(0,0)[lb]{$\vep$}}%
\put(23.1000,-21.9000){\makebox(0,0)[lb]{$\vep/2$}}%
%
{\color[named]{Black}{%
\special{pn 8}%
\special{pa 3150 2590}%
\special{pa 3150 610}%
\special{dt 0.045}%
}}%
%
{\color[named]{Black}{%
\special{pn 8}%
\special{pa 3710 2600}%
\special{pa 3710 620}%
\special{dt 0.045}%
}}%
%
{\color[named]{Black}{%
\special{pn 8}%
\special{pa 2590 2150}%
\special{pa 4880 2150}%
\special{dt 0.045}%
}}%
%
{\color[named]{Black}{%
\special{pn 8}%
\special{pa 2610 1640}%
\special{pa 4900 1640}%
\special{dt 0.045}%
}}%
\put(40.3000,-25.1000){\makebox(0,0)[lb]{$R_{1,\infty}\times R_{2,0}$}}%
\put(26.9000,-11.5000){\makebox(0,0)[lb]{$R_{1,0}\times R_{2,\infty}$}}%
\put(40.3000,-11.3000){\makebox(0,0)[lb]{$R_{1,\infty}\times R_{2,\infty}$}}%
%
{\color[named]{Black}{%
\special{pn 4}%
\special{pa 3560 2150}%
\special{pa 2610 1200}%
\special{fp}%
\special{pa 3680 2150}%
\special{pa 2610 1080}%
\special{fp}%
\special{pa 3710 2060}%
\special{pa 2610 960}%
\special{fp}%
\special{pa 3710 1940}%
\special{pa 2610 840}%
\special{fp}%
\special{pa 3710 1820}%
\special{pa 2610 720}%
\special{fp}%
\special{pa 3710 1700}%
\special{pa 2620 610}%
\special{fp}%
\special{pa 3710 1580}%
\special{pa 2730 600}%
\special{fp}%
\special{pa 3710 1460}%
\special{pa 2850 600}%
\special{fp}%
\special{pa 3710 1340}%
\special{pa 2970 600}%
\special{fp}%
\special{pa 3710 1220}%
\special{pa 3090 600}%
\special{fp}%
\special{pa 3710 1100}%
\special{pa 3210 600}%
\special{fp}%
\special{pa 3710 980}%
\special{pa 3330 600}%
\special{fp}%
\special{pa 3710 860}%
\special{pa 3450 600}%
\special{fp}%
\special{pa 3710 740}%
\special{pa 3570 600}%
\special{fp}%
\special{pa 3440 2150}%
\special{pa 2610 1320}%
\special{fp}%
\special{pa 3320 2150}%
\special{pa 2610 1440}%
\special{fp}%
\special{pa 3200 2150}%
\special{pa 2610 1560}%
\special{fp}%
\special{pa 3080 2150}%
\special{pa 2610 1680}%
\special{fp}%
\special{pa 2960 2150}%
\special{pa 2610 1800}%
\special{fp}%
\special{pa 2840 2150}%
\special{pa 2610 1920}%
\special{fp}%
\special{pa 2720 2150}%
\special{pa 2610 2040}%
\special{fp}%
}}%
%
{\color[named]{Black}{%
\special{pn 4}%
\special{pa 4230 1640}%
\special{pa 3270 2600}%
\special{fp}%
\special{pa 4110 1640}%
\special{pa 3160 2590}%
\special{fp}%
\special{pa 3990 1640}%
\special{pa 3150 2480}%
\special{fp}%
\special{pa 3870 1640}%
\special{pa 3150 2360}%
\special{fp}%
\special{pa 3750 1640}%
\special{pa 3150 2240}%
\special{fp}%
\special{pa 3630 1640}%
\special{pa 3150 2120}%
\special{fp}%
\special{pa 3510 1640}%
\special{pa 3150 2000}%
\special{fp}%
\special{pa 3390 1640}%
\special{pa 3150 1880}%
\special{fp}%
\special{pa 3270 1640}%
\special{pa 3150 1760}%
\special{fp}%
\special{pa 4350 1640}%
\special{pa 3390 2600}%
\special{fp}%
\special{pa 4470 1640}%
\special{pa 3510 2600}%
\special{fp}%
\special{pa 4590 1640}%
\special{pa 3630 2600}%
\special{fp}%
\special{pa 4710 1640}%
\special{pa 3750 2600}%
\special{fp}%
\special{pa 4830 1640}%
\special{pa 3870 2600}%
\special{fp}%
\special{pa 4940 1650}%
\special{pa 3990 2600}%
\special{fp}%
\special{pa 4950 1760}%
\special{pa 4110 2600}%
\special{fp}%
\special{pa 4950 1880}%
\special{pa 4230 2600}%
\special{fp}%
\special{pa 4950 2000}%
\special{pa 4350 2600}%
\special{fp}%
\special{pa 4950 2120}%
\special{pa 4470 2600}%
\special{fp}%
\special{pa 4950 2240}%
\special{pa 4590 2600}%
\special{fp}%
\special{pa 4950 2360}%
\special{pa 4710 2600}%
\special{fp}%
\special{pa 4950 2480}%
\special{pa 4830 2600}%
\special{fp}%
}}%
%
{\color[named]{Black}{%
\special{pn 4}%
\special{pa 4980 1590}%
\special{pa 3140 1590}%
\special{fp}%
\special{pa 4980 1710}%
\special{pa 3140 1710}%
\special{fp}%
\special{pa 4980 1830}%
\special{pa 3140 1830}%
\special{fp}%
\special{pa 4980 1950}%
\special{pa 3140 1950}%
\special{fp}%
\special{pa 4980 2070}%
\special{pa 3140 2070}%
\special{fp}%
\special{pa 4980 1470}%
\special{pa 3140 1470}%
\special{fp}%
\special{pa 4980 1350}%
\special{pa 3140 1350}%
\special{fp}%
\special{pa 4980 1230}%
\special{pa 3140 1230}%
\special{fp}%
\special{pa 4980 1110}%
\special{pa 3140 1110}%
\special{fp}%
\special{pa 4980 990}%
\special{pa 3140 990}%
\special{fp}%
\special{pa 4980 870}%
\special{pa 3140 870}%
\special{fp}%
\special{pa 4980 750}%
\special{pa 3140 750}%
\special{fp}%
\special{pa 4980 630}%
\special{pa 3140 630}%
\special{fp}%
}}%
\end{picture}}%
\caption{Open covering $R_{\infty}^{\rm prod}$.}
\label{proddomein}
\end{figure}
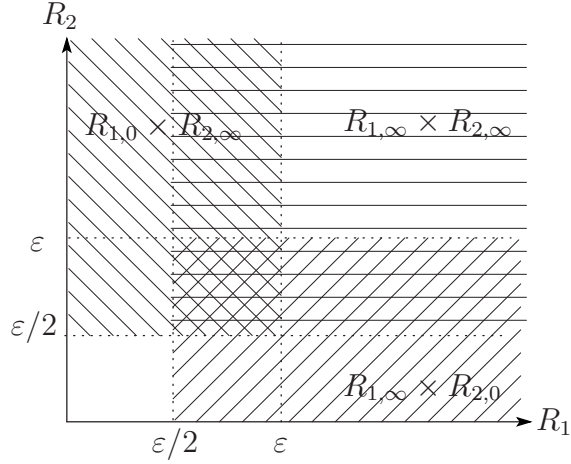

\begin{prop}\label{cobandprod}
We have the following equality among equivariant local indices. 
$$
\ind_{G_1\times G_2}(R, \tilde R_{\infty})=\ind_{G_1\times G_2}(R, R_{\infty}^{\rm prod})=
\ind_{G_1}(R_1,R_1^{\times})\otimes\ind_{G_2}(R_2,R_2^{\times})\in R(G_1\times G_2).  
$$
\end{prop}
\begin{proof}
The first equality follows from the cobordism invariance of local index (\cite[Theorem~7.1]{Fujitacobinv}). In fact the union of these two acyclic compatible systems on $\tilde R_{\infty}$ and $R_{\infty}^{\rm prod}$ is also $G_1\times G_2$-tangential acyclic compatible system. The second equality follows from the product formula (\cite[Theorem~8.8]{Fujita-Furuta-Yoshida2}). 
\end{proof}
\section{Local index of folded cylinder}
\label{A computation of local index of the folded cylinder}

In this appendix we consider a natural folded symplectic structure on the cylinder and several geometric structures on it, which plays important role in the study of local property of the neighbourhood of the fold in a folded symplectic manifold. We consider a perturbation of the Dirac operator and give the direct computation of the $L^2$-kernel of the perturbed Dirac operator. We show that the $L^2$-kernel is trivial, in particular, the local index is equal to $0$. 

For any $\vep >0$, a folded symplectic structure on a cylinder (of finite length) $M_{\vep}:=(-\vep,\vep)\times S^1$ is given by a closed 2-form $2rdr\wedge d\theta$, where $(r,\theta)$ is a coordinate function on $M_{\vep}$. Here we use the opposite orientation of the cylinder as that in Section~\ref{Compatible system on toric origami manifolds} and subsequent argument for conventional reason. The standard $S^1$-action on the $S^1$-factor is Hamiltonian (in fact it is toric origami) with the moment map $(r,\theta)\mapsto r^2$. Moreover the trivial line bundle $L_0$ with connection $d-2\pi\sqrt{-1}r^2d\theta$ and the trivial lift of the $S^1$-action to the fiber direction gives an $S^1$-equivariant pre-quantizing line bundle over $M_{\vep}$. To give a computation of the local index of this toric origami manifold, we need a Clifford module bundle, Dirac-type operator along the $S^1$-orbits over a completion of $M_{\vep}$ as a Riemannian manifold.  
We summarize the set-up as follows. 

\medskip

{\bf SET-UP.}
\begin{itemize}
\item $M:=\R\times S^1$ : cylinder of infinite length
\item $(r,\theta)$ : coordinate function on $M$
\item $g:=dr^2+d\theta^2$ : Riemannian metric on $M$
\item $\rho : \R\to \R$ : smooth function with 
$$
\rho(r)=\left\{ 
\begin{array}{lll}
r^2 \quad (|r| < 1/4)) \\ 
 1/2  \quad (|r|>1/2) 
\end{array}\right.
$$
\item $\omega:=\rho'(r)dr\wedge d\theta$ : closed 2-form on $M$
\item $J:\partial_r\mapsto\partial_{\theta}, \ \partial_{\theta}\mapsto -\partial_r$ : almost complex structure on $M$
\item $TM_{\C}=(TM,J)$ : complex tangent bundle with frame $\partial_{\theta}$
\item $W^+:=M\times\C$, $W^-:=TM_{\C}$, $W:=W^{+}\oplus W^{-}$ : $\Z/2$-graded vector bundle
\item $c:T^*M\to {\rm End}(W)$ : Clifford action on $W$ defined by 
$$c(dr)=
\begin{pmatrix}
0 & -\sqrt{-1} \\ 
-\sqrt{-1} & 0 
\end{pmatrix}, \quad 
c(d\theta)=
\begin{pmatrix}
0 & -1 \\ 
1 & 0
\end{pmatrix} 
$$
\item $\nabla^W=d-2\pi\rho(r)
\begin{pmatrix}
1 & 0 \\ 
0 & 1
\end{pmatrix}d\theta$ : 
Clifford connection of $W$ 
\item 
$D=D^++D^-:\Gamma(W)\to \Gamma(W)$ : Dirac operator, 
\begin{eqnarray*}
D&=&c(\partial_r)\nabla^W_{\partial_r}+c(\partial_{\theta})\nabla^W_{\partial_{\theta}}=D^++D^- \\ 
&=&\begin{pmatrix}
0 &  -\partial_{\theta}-\sqrt{-1}\partial_r+2\pi\sqrt{-1}\rho \\
\partial_{\theta}-\sqrt{-1}\partial_r-2\pi\sqrt{-1}\rho & 0
\end{pmatrix}
\end{eqnarray*}
\item Let $S^1$ acts on $M$ in the standard way, and   
we take a lift of 
the $S^1$-action on $W$ so that the action on the fiber direction is trivial. 

\item $D_{S^1}=D_{S^1}^++D_{S^1}^-:\Gamma(W) \to \Gamma(W)$ : Dirac operator along the $S^1$-orbits,  
$$
D_{S^1}=c(\partial_{\theta})\nabla_{\partial_{\theta}}^W=
\begin{pmatrix}
0 &  -\partial_{\theta}+2\pi\sqrt{-1}\rho \\
\partial_{\theta}-2\pi\sqrt{-1}\rho & 0
\end{pmatrix}
$$
\end{itemize}


\begin{remark}
When we consider the restriction to an $S^1$-invariant small open neighbourhood  $M_{\vep}$ of $\{0\}\times S^1=S^1$ in $M$, the closed 2-form $\omega$ is the folded symplectic form on $M_{\vep}$ and the $\Z/2$-graded Clifford module bundle $W$ is the one associated with the pre-quantizing line bundle $L_0$, the trivial line bundle with the connection $d-2\pi\sqrt{-1}\rho d\theta$. Note that $M_{\vep}$ has a unique spin$^c$-structure and the Clifford module bundle $W$ which is isomorphic to $W_{0,L_0}:={\rm Hom}_{Cl_2}(W_2, \wedge_{\C}^{\bullet}(TS^1\oplus\R^3))\otimes L_0$ as in Proposition~\ref{product structure of Cl}.  
\end{remark}

By using this data we have a compatible system on $M_{\vep}$ and can define the $S^1$-equivariant local index $\ind_{S^1}(S^1\times(-\vep,\vep), S^1\times(-\vep,\vep)\setminus S^1)$. The index is defined by the following perturbation of the Dirac operator:  
$$
D_t=D_t^++D_t^-,  \ D_t^+:=D^++tD_{S^1}^+,  \ D_t^-:=D^-+tD_{S^1}^-, 
$$
$$
D_t^+=(1+t)(\partial_{\theta}-2\pi\sqrt{-1}\rho)-\sqrt{-1}\partial_r, 
$$and 
$$
D_t^-=-(1+t)(\partial_{\theta}-2\pi\sqrt{-1}\rho)-\sqrt{-1}\partial_r.   
$$

\begin{prop}\label{vanishingoffoldedcylinder}
We have $\ker_{L^2}(D_t^+)=\ker_{L^2}(D_t^-)=0$ for any $t\geq 0$. In particular we have  $\ind_{S^1}(S^1\times(-\vep,\vep), S^1\times(-\vep,\vep)\setminus S^1)=0$, for any $\vep>0$. 
\end{prop}

\begin{proof}
By using the Fourier expansion $\displaystyle\phi(r,\theta)=\sum_{m\in\Z}a_m(r)e^{2\pi\sqrt{-1}m\theta}$ for smooth section $\phi$ of $W^+$, the equation $D_t^+\phi=0$ can be rewritten as a series of differential equations 
$$
a'_m(r)=2\pi(1+t)(m-\rho(r))a_m(r) \quad (m\in\Z).
$$ Each of these equations has solutions 
$$
a_m(r)=\alpha_m\exp\left(2\pi(1+t)\int_0^r(m-\rho(r))dr\right),
$$ where $\alpha_m\in\C$ is constant. 
Suppose that the solution $\phi$ is an $L^2$-section. Since $\rho\equiv 1/2$ on $\pm r\gg 0$ we have $\alpha_m=0$ for all $m\in\Z$. 
In particular there are no non-trivial $L^2$-solutions of $D_t^+\phi=0$ for any $t$. As in the same way the equation $D_t^-\phi=0$ for $\displaystyle\phi(r,\theta)=\sum_{m\in\Z}b_m(r)e^{2\pi\sqrt{-1}m\theta}$ has solutions 
$$
b_m(r)=\beta_m\exp\left(-2\pi(1+t)\int_0^r(m-\rho(r))dr\right) \quad (m\in\Z), 
$$ for any constant $\beta_m\in\C$ and one can see that there are no non-trivial $L^2$-solutions. 
\end{proof}

\begin{remark}
The vanishing of the index can be deduced from the existence of 
an orientation reversing isomorphism of $S^1\times (-\vep, \vep)$ defined by $(\theta,t)\mapsto (\theta, -t)$.  
\end{remark}


\vspace{1cm}

\noindent{\bf Acknowledgements.}
The author would like to thank Mikio Furuta and Takahiko Yoshida for stimulating  conversations. Especially the argument in Section~\ref{Toric case} is based on the discussion with them. 
The auothor is grateful to the anonymous referee for his/her comments and pointing out some mistakes on the previous version. Due to referee's comments on the proof of Theorem~\ref{crackind=0}, the author could simplify it. 

\bibliographystyle{amsplain}
\bibliography{reference}

\end{document}